\documentclass{m2an}
\usepackage{graphicx}
\newtheorem{rmrks}[thrm]{Remarks}
\numberwithin{equation}{section}
\begin{document}
\title{Influence of Bottom Topography on Long Water Waves}

\author{Florent Chazel}\address{Laboratoire de Math\'ematiques
  Appliqu\'ees de Bordeaux, Universit\'e Bordeaux 1, 351 Cours de la Lib\'eration,
F-33405 Talence cedex ; \email{florent.chazel@math.u-bordeaux1.fr}}
\date{}
\begin{abstract} 
We focus here on the water
waves problem for uneven bottoms in the long-wave regime, on an
        unbounded two or three-dimensional domain. In order to derive
        asymptotic models for this problem, we consider two different regimes of bottom
        topography, one for small variations in amplitude, and one for
        strong variations. Starting from the Zakharov
        formulation of this problem, 
        we rigorously compute the asymptotic
        expansion of the involved Dirichlet-Neumann operator. Then,
        following the global strategy introduced by Bona, Colin and
        Lannes in \cite{BCL}, new symetric asymptotic models are derived for each
        regime of bottom topography. Solutions of these systems are
        proved to give good approximations of solutions of the water
        waves problem. These results hold for solutions that
        evanesce at infinity as well as for spatially periodic ones. \end{abstract}
\begin{resume}
Nous nous int\'eressons ici au probl\`eme d'Euler surface
libre pour des fonds non plats en r\'egime d'ondes longues, sur un
domaine non born\'e \`a deux ou trois dimensions. Afin de
construire des mod\`eles asymptotiques pour ce probl\`eme, nous consid\`erons deux r\'egimes
topographiques sur le fond du domaine, l'un pour de petites variations 
en amplitude, et l'autre pour de fortes variations. A
partir de la formulation de Zakhzarov, nous
contruisons rigoureusement le
d\'eveloppement asymptotique de l'op\'erateur de Dirichlet-Neumann relatif
au probl\`eme. En suivant la strat\'egie globale propos\'ee
par Bona, Colin et Lannes dans \cite{BCL}, nous obtenons ensuite de nouveaux
mod\`eles asymptotiques sym\'etriques pour chaque r\'egime de variation
topographique du fond. Nous prouvons alors que les solutions de ces
syst\`emes fournissent de bonnes approximations aux solutions des
\'equations d'Euler surface libre. Ces r\'esultats sont valables aussi
bien pour des solutions \'evanescentes \`a l'infini que pour des solutions
spatialement p\'eriodiques. \end{resume}

\maketitle
\section*{Introduction}
	
        \vspace{1em}
	
	\subsection*{Generalities}
	
	This paper deals with the water waves problem for uneven bottoms
	which consists in describing the motion of
	the free surface and the evolution of the 
	velocity field of a layer of fluid, under the following assumptions : the fluid is ideal, 
	incompressible, irrotationnal, and under the only 
	influence of gravity. \\ Earlier works have set a good theoretical background for this problem : its 
	well-posedness has been discussed among others by Nalimov (\cite{Nalimov}, 1974), Yoshihara 
	(\cite{Yoshihara1}, 1982), Craig (\cite{Craig}, 1985), Wu (\cite{Wu1}, 1997, \cite{Wu2}, 1999) and 
	Lannes (\cite{Lannes}, 2005).\\
	Nevertheless, the solutions of these equations are very
        difficult to describe, because of the complexity of these
        equations. At this point, a classical method is to choose an
        asymptotic regime, in which we look for approximate
        models and hence for approximate solutions.
	We consider in this paper the so-called long-wave regime,
        where the ratio of the typical 
	amplitude of the waves over the mean depth and the ratio of
        the square of the mean depth over the square of the typical wave-length are both
        neglictible in front of 1 and of the same order.\\
        In 2002, Bona, Chen and Saut constructed in \cite{BCS} a
	large class of systems for this regime and performed a formal study in the
        two-dimensional case.
	A significant step forward has been made in 2005 by Bona,
        Colin and Lannes in \cite{BCL}; they rigorously justified
	the systems of Bona, Chen and Saut, and derived a new specific
        class of symmetric systems. Solutions of these systems are
        proved to
	tend
	to solutions of the water waves problem on a long time scale,
        as the amplitude becomes small and the wavelength large. Thanks to their symmetric structure, 
	computing solutions 
	of such systems is significantly easier than computing directly solutions of the water waves
	problem. Another significant work in this field is the one of
        Lannes and Saut (\cite{LannesJCS}, 2006) on weakly transverse
        Boussinesq systems. 
        \\ However, all these results only hold for flat
        bottoms. The case of uneven bottoms has been less investigated ;
        some of the significant references are Peregrine (\cite{Peregrine}, 1967), Madsen et
        al. (\cite{Madsen}, 1991), Nwogu (\cite{Nwogu}, 1993), and Chen
        (\cite{Chen}, 2004). Peregrine was the first one to formulate
        the classical Boussinesq equations for waves in shallow water
        with variable depth on a three-dimensionnal domain. Following
        this work, Madsen et al. and Nwogu derived new Boussinesq-like
        systems for uneven bottoms with improved linear dispersion
        properties. Recently, Chen performed a formal
        study of the water waves problem for uneven bottoms with small
        variations in amplitude, in 1D of
        surface, and derived a class of asymptotic models
        inspired by the work of Bona, Chen and Saut. 
        To our knowledge, the only rigorously justified result on the uneven bottoms case
        is the work of Iguchi (\cite{Iguchi}, 2004), who provided a rigorous approximation via a system of KdV-like 
	equations, in the case of a slowly varying bottom.\\ 
        The main idea of our paper is to reconsider
	the water waves problem for uneven bottoms in the angle shown
        by Bona, Colin and Lannes. Moreover, our goal is to consider two different types of
        bottoms : bottoms with small variations in amplitude, and
        bottoms with strong variations in amplitude.  
	To this end, we introduce a new parameter to characterize the
        shape of the bottom. In the end, new asymptotic models are
        derived, studied and rigorously justified under the assumption
        that long time solutions to the water waves equations exist.
	
        \vspace{1em}
	
	\subsection*{Presentation and formulation of the problem}
	
	In this paper, we work indifferently in two or three dimensions. 
	Let us denote by $X\in\xR^d$ the transverse
	variable, $d$ being equal to 1 or 2. 
	In the two-dimensional case, $d=1$ and $X=x$ corresponds to the coordinate
	along the primary direction of propagation whilst in the three-dimensional case, 
	$d=2$ and $X=(x,y)$ represents
	the horizontal variables. We restrict our study to the case where the free surface and the bottom
	can be described by the graph
	of two functions $(t,X) \rightarrow \eta(t,X)$ and $X \rightarrow b(X)$ defined respectively 
	over the surface $z=0$ and the mean depth $z=-h_0$ both at the steady state, 
	$t$ corresponding to the time variable. 
	The time-dependant domain $\Omega_t$ of the fluid is thus taken of the form : 
	$$\Omega_t = \{(X,z),\,X\in\xR^d,\,-h_0+b(X) \le z\le\eta(t,X)\}\;\;.$$
	
	\vspace{1em}
	
        \begin{figure}
        \centering \includegraphics[height=0.4\hsize]{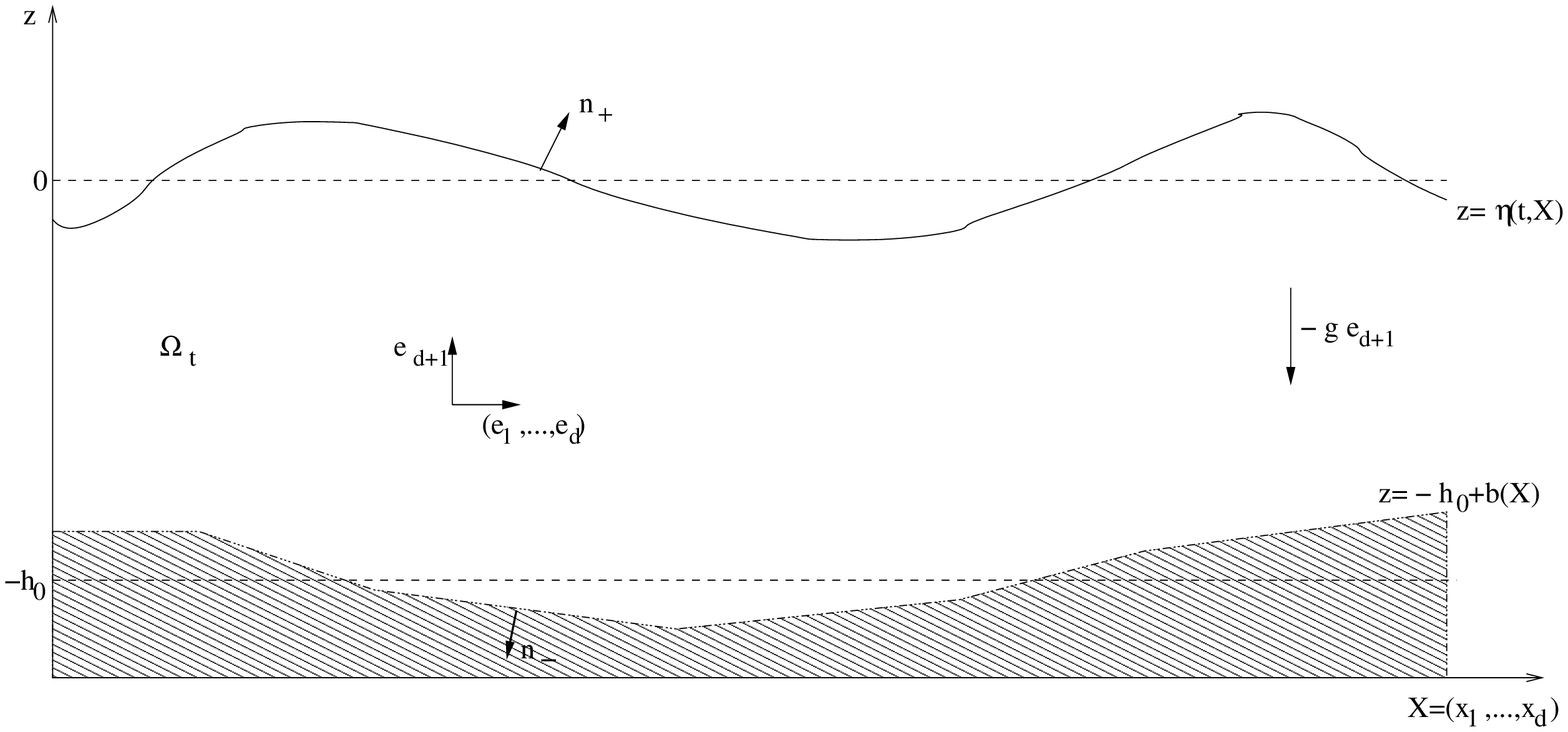}
        \end{figure}
	
	\vspace{1em}
	
	\noindent In order to avoid some special
	physical cases
	such as the presence of islands or beaches, we set a condition of minimal water depth : there exists
	a strictly positive constant $h_{min}$ such that
	\begin{equation} \label{hmin}
	\eta(t,X) + h_0 - b(X)  \ge h_{min}\;\;,\;\; (t,X) \in \xR\times\xR^2\;\;.
	\end{equation}
	For the sake of simplicity, we assume here that $b$ and all
        its derivatives are bounded.
		
	\vspace{1em}
	
	\noindent The motion of the fluid is described by the following system of equations :\vspace{1em}
	\begin{equation} \label{V}
	\left\{
	\begin{array}{cl}
	\vspace{1em}
	\partial_t V + V\,\cdot\nabla_{X,z} V = -g e_z - \nabla_{X,z} P & \;\;\;\mbox{in}\; \Omega_t, \;\; t \ge 0\;\;, \\
	\vspace{1em}
	\nabla_{X,z} \cdot V = 0 & \;\;\;\mbox{in}\; \Omega_t, \;\; t \ge 0\;\;, \\
	\vspace{1em}
	\nabla_{X,z} \times V = 0 & \;\;\;\mbox{in}\; \Omega_t, \;\; t \ge 0\;\;, \\
	\vspace{1em}
	\partial_t \eta	- \sqrt{1+|\nabla_X \eta|^2}\;\mbox{{\bfseries n}}_+ \cdot V|_{z=\eta (t,X)} = 0 & \;\;\; \mbox{for}\;
	t \ge 0,\;X\in\xR^d\;\;,\\ 
	\vspace{1em}
	P_{|_{z=\eta(t,X)}} = 0 & \;\;\; \mbox{for}\;t \ge 0,\;X\in\xR^d\;\;,\\
	\vspace{1em}
	\mbox{{\bfseries n}}_- \cdot V|_{z=-h_0+b(X)} = 0 & \;\;\; \mbox{for}\; t \ge 0,\;X\in\xR^d\;\;,\\
	\end{array}
	\right.
	\end{equation}
	\vspace{1em}
	\noindent where $\mathbf{n}_+ = \frac{1}{\sqrt{1+|\nabla
	\eta|^2}} (-\nabla \eta, 1)^T$
	denotes the outward normal vector to the surface and $\mathbf{n}_- = 	 
	\frac{1}{\sqrt{1+|\nabla
	b|^2}} (\nabla b, -1)^T$ denotes the outward normal vector to the bottom. 	
	The first equation corresponds to the Euler equation for a perfect fluid under the influence 
	of gravity
	(which is characterized by the term $-g e_z$ where $e_z$ denotes the base vector 
	along the vertical component). The
	second and third one characterize the incompressibility and irrotationnality 
	of the fluid. The
	fourth and last ones deal with the boundary conditions at the surface and the bottom. 
	These are given by the usual assumption
	that they are both bounding surfaces, i.e. surfaces across which no fluid particles are transported.
	As far as the pressure $P$ is concerned, we assume that it is constant at the 
	surface by neglicting the surface tension. Up to a renormalization, we can assume that it is equal to 
	zero at the surface.
	
	\vspace{1em}
	
	In this paper, we use the Bernoulli formulation of the water-waves equations. 
	The conditions of incompressibility
	and irrotationnality ensure the existence of a potential flow $\phi$ such that $V=\nabla_{X,z}\phi$. 
	From now on, we
	separate the transverse variable $X \in \xR^d$ and the vertical variable $z \in \xR$ : the
	operators $\nabla$ and $\Delta$  act only on the transverse variable $X \in \xR^d$ so that we have 
	$V=\nabla\phi + \partial_z^2 \phi$.
	The use of the potential flow $\phi$ instead of the velocity $V$ leads to the following formulation of 
	(\ref{V}) :
	\vspace{1em}
	\begin{equation} \label{F}
	\left\{
	\begin{array}{cl}
	\vspace{1em}
	\partial_t \phi + \frac{1}{2}\,\left[\,| \nabla \phi |^2 + | \partial_z \phi |^2 \,\right] + g z = -P & 
	\;\;\;\mbox{in}\; \Omega_t, \;\; t \ge 0\;\;, \\
	\vspace{1em}
	\Delta \phi + \partial_z^2 \phi = 0 & \;\;\;\mbox{in}\; \Omega_t, \;\; t \ge 0\;\;, \\
	\vspace{1em}
	\partial_t \eta	- \sqrt{1+|\nabla \eta|^2}\;\partial_{\mbox{{\bfseries n}}_+} \phi |_{z=\eta (t,X)} = 0 & \;\;\; 
	\mbox{for}\; t \ge 0,\;X\in\xR^d\;\;,\\ 
	\vspace{1em}
	\partial_{\mbox{{\bfseries n}}_-} \phi |_{z=-h_0+b(X)} = 0 & \;\;\; \mbox{for}\; t \ge 0,\;X\in\xR^d\;\;,\\
	\end{array}
	\right.
	\end{equation}
	
	\noindent where we used the notations 
	$\partial_{\mbox{{\bfseries n}}_-} = \mbox{{\bfseries n}}_- \cdot 
	\left( \begin{array}{l}
	\nabla \\ 
	\partial_z
	\end{array}\right)$ and $\partial_{\mbox{{\bfseries n}}_+} = \mbox{{\bfseries n}}_+ \cdot \left( \begin{array}{l}
	\nabla \\ 
	\partial_z
	\end{array}\right)$.
	
	\vspace{1em}
	\noindent Separating the 
	variables $X$ and $z$ in the boundary conditions and taking the trace
	of (\ref{F}) on the free surface thus leads to the system :
	\vspace{1em}
	\begin{equation} \label{G}
	\left\{
	\begin{array}{cl}
	\vspace{1em}
	\Delta \phi + \partial_z^2 \phi = 0 & \;\;\;\mbox{in}\; \Omega_t, \;\; t \ge 0\;\;, \\
	\vspace{1em}
	\partial_t \phi + \frac{1}{2}\,\left[\,| \nabla \phi |^2 + | \partial_z \phi |^2 \,\right] + g \eta = 0 & 
	\;\;\;\mbox{at}\; z = \eta(t,X),\;X \in \xR^d,\; t \ge 0\;\;, \\
	\vspace{1em}
	\partial_t \eta	+ \nabla \eta \cdot \nabla \phi - \partial_z \phi = 0 & \;\;\; 
	\mbox{at}\; z = \eta(t,X),\;X \in \xR^d,\; t \ge 0\;\;, \\
	\vspace{1em}
	\nabla b \cdot \nabla \phi - \partial_z \phi = 0 & \;\;\; \mbox{at}\; z = -h_0+b(X),\;X \in \xR^d,\; t \ge 0\;\;. \\
	\end{array}
	\right.
	\end{equation}
	
	We now perform a non-dimensionalisation of these equations using the following 
	parameters : $\lambda$ is the typical
	wavelength, $a$ the typical amplitude of the waves, $h_0$ the mean depth of the fluid, $b_0$ the
	typical amplitude of the bottom, 
	$t_0=\frac{\lambda}{\sqrt{g h_0}}$ a typical period of time ($\sqrt{g h_0}$ corresponding
	to sound velocity in the fluid) and $\phi_0 = \frac{\lambda a}{h_0}\,\sqrt{g h_0}$.
	Introducing the following parameters :
	$$ \epsilon = \frac{a}{h_0} \;\;;\;\; \beta=\frac{b_0}{h_0} \;\;;\;\; S=\frac{a \lambda^2}{h_0^3}\;\;, $$ 
	and taking the Stokes number $S$ to be equal to one, 
	one gets for the non-dimensionnalized version of (\ref{G}) :
	\vspace{1em}
	\begin{equation} \label{H}
	\left\{
	\begin{array}{cl}
	\vspace{1em}
	\varepsilon\,\Delta \phi + \partial_z^2 \phi = 0 & \;\;\;-1+\beta\, b \le z 
	\le \varepsilon \eta , \;
	X \in \xR^d , \; t \ge 0\;\;, \\
	\vspace{1em}
	\partial_t \phi + \frac{1}{2}\,\left[\,\varepsilon | \nabla \phi |^2 + | \partial_z 
	\phi |^2 \,\right] + g \eta = 0
	& \;\;\;\mbox{at}\; z = \varepsilon\eta,\;X \in \xR^d,\; t \ge 0\;\;, \\
	\vspace{1em}
	\partial_t \eta	+ \varepsilon \,\nabla \eta \cdot \nabla \phi - \frac{1}{\varepsilon} 
	\partial_z \phi = 0 & \;\;\; 
	\mbox{at}\; z = \varepsilon\eta,\;X \in \xR^d,\; t \ge 0\;\;, \\
	\vspace{1em}
	\partial_z \phi - 	\varepsilon \beta \, \nabla b \cdot \nabla \phi = 0 
	& \;\;\; \mbox{at}\; z = 
	-1+\beta\, b,	
	\;X \in \xR^d,\; t \ge 0\;\;. \\
	\end{array}
	\right.
	\end{equation}
	
	The final step consists in recovering the Zakharov formulation by
        reducing the previous system (\ref{H}) to a system expressed
	at the free surface. To this end, we introduce the trace of the velocity potential $\phi$ 
	at the free surface, namely $\psi$ :
	$$ \psi(t,X) = \phi(t,X,\varepsilon\,\eta(t,X))\;\;,$$
	and the operator $Z_{\varepsilon}(\varepsilon\eta,\beta b)$ which maps $\psi$ to $\partial_z \phi 
	|_{z=\varepsilon\,\eta}$. This operator is defined
	for any $f \in (C^1 \cap W^{1,\infty})(\xR^d)$ by :
	\begin{equation} \label{Ze}
	Z_{\varepsilon}(\varepsilon\eta,\beta b)f :
	\left(
	\begin{array}{cll}
	\vspace{0.5em}
	H^{\frac{3}{2}}(\xR^d) & \longrightarrow & H^{\frac{1}{2}}(\xR^d)\\
	\vspace{0.5em}
	f & \longmapsto & \partial_z u_{|_{z=\varepsilon\eta}} \;\;\;\mbox{with } u \mbox{ solution of :}\\
	\vspace{0.5em}
	& & \varepsilon\,\Delta u + \partial_z^2 u = 0, \;\;\;-1+\beta\, b \le z \le \varepsilon \eta \;,\\
	\vspace{0.5em}
	& & \partial_z u - \varepsilon \beta\,\nabla b \cdot \nabla u = 0,\;\;\;z=-1+\beta b\;,\\
	\vspace{0.5em}
	& & u(X,\varepsilon \eta) = f, \;\;\; X \in \xR^d\;.\\
	\end{array}
	\right)\;\;.
	\end{equation}
	
	\noindent Using this operator and computing the derivatives of $\psi$ in terms of $\psi$ and $\eta$,
	the final formulation $(S_{0})$ of the water waves problem reads :
	\begin{equation} \label{S0}
	(S_{0})\left\{
	\begin{array}{l}
	\vspace{1em}
	\partial_t \psi - \varepsilon \partial_t \eta Z_{\varepsilon}(\varepsilon\eta,\beta b)\psi
	+\frac{1}{2}\,\Big[\,\varepsilon
	\,|\nabla \psi - \varepsilon\,\nabla \eta
        Z_{\varepsilon}(\varepsilon\eta,\beta b)\psi |^2 
        + \,|Z_{\varepsilon}(\varepsilon\eta,\beta b)\psi |^2
	\, \Big] + \eta = 0 \;\;,\\
	\partial_t \eta + \varepsilon \nabla \eta \cdot \left[\,\nabla \psi - \varepsilon \nabla 
	\eta Z_{\varepsilon}(\varepsilon\eta,\beta b)
	\psi \,\right] = \frac{1}{\varepsilon}\,Z_{\varepsilon}(\varepsilon\eta,\beta b)\psi\;\;.
	\end{array}
	\right.
	\end{equation}
	
        \vspace{1em}

	\subsection*{Organization of the paper}
	
	The aim of this paper is to derive and study two different asymptotic regimes based 
	each on a specific assumption
	on the parameter $\beta$ which characterizes the topography of the bottom. 
	The first assumption deals with the case
	$\beta = O(\varepsilon)$ which corresponds to the physical case of a bottom with 
	small variations in amplitude. The
	second one deals with the more complex case $\beta = O(1)$ which corresponds 
	to the physical case of a bottom
	with high variations in amplitude. \\
	The following part will be devoted to the asymptotic expansion of the operator 
	$Z_{\varepsilon}(\varepsilon\eta,\beta b)$ in the two regimes
	mentionned above. To this end, a general method is introduced and rigorously proved
	which aims at deriving asymptotic expansions of Dirichlet-Neumann operators for a large
	class of elliptic problems.
	This result is then applied in each regime, wherein a formal expansion is performed and an asymptotic 
	Boussinesq-like model of (\ref{S0}) is derived.
	The second and third part are both devoted to the derivation of new classes of equivalent systems, following
	the strategy developped in \cite{BCL}. In the end, completely symmetric systems are obtained for
	each bottom topography regime : convergence results are proved showing that solutions of these symmetric asymptotic
	systems tend to associated solutions of the water waves problem.
	
        \vspace{2em}

	\section{Asymptotic expansion of the operator $Z_{\varepsilon}(\varepsilon\eta,\beta b)$}
	
        \vspace{1em}
	
	\noindent This section is devoted to the asymptotic expansion
        of the operator $Z_\varepsilon(\varepsilon \eta,\beta b)$ defined in the
        previous section as $\varepsilon$ tends to zero, in
        both regimes $\beta = O(\varepsilon)$ and $\beta = O(1)$. 
		To this end, we first enounce some general results on elliptic
        equations on a strip : the final proposition gives a general rigourously justified method for
        determining an approximation of Dirichlet-Neumann
        operators. This result is then applied to
        the case of the operator $Z_{\varepsilon}(\varepsilon\eta,\beta b)$ and two asymptotic
        models with bottom effects are derived.
	
	\subsection{Elliptic equations on a strip}
	
	\noindent In this part, we aim at studying a general elliptic equation on
	a domain $\Omega$ given by :
	$$ \Omega = \{(X,z) \in \xR^{d+1} / X \in \xR^d, -h_0+B(X) < z < \eta(X)\}\;\;, $$
	where the functions $B$ and $\eta$ satisfy the following condition :
	
	\begin{equation} \label{cond} 
        \exists \, h_{min} > 0 \;,\;\forall X \in \xR^d \;,\; \eta(X)-B(X)+h_0 \ge h_{min}  \;\;.   
        \end{equation}

        \vspace{1em}
    
	\noindent Let us consider the following general elliptic boundary
        value problem set on the domain $\Omega$ :
	
	\begin{equation} \label{General}	
	-\nabla_{X,z}\,.\,P\,\nabla_{X,z}\,u = 0 \;\;\;\;\; \mbox{in}\; \Omega	\;\;,
        \end{equation}
	\begin{equation}	 \label{bord}	
	u_{\,|_{z=\eta(X)}} = f \;\;\; \mbox{and} \;\;\; \partial_n\,
        u_{\,|_{z=
	-h_0+B(X)}} = 0\;\;,
	\end{equation}

        \vspace{1em}
    	
	\noindent where $P$ is a diagonal $(d+1)\times(d+1)$ matrix
        whose coefficients $(p_i)_{1\le i\ \le d+1}$ are 
	constant and 
	strictly positive. Straightforwardly $P$ is 
        coercive. We denote by 
        $\partial_n\,u_{\,|z=-h_0+B(X)}$ the outward conormal derivative associated to 
        $P$ of $u$
	at the lower boundary $\{z=-h_0+B(X)\}$, namely :
	$$ \partial_n\,u_{\,|_{z=-h_0+B(X)}} = 	
	-\mathbf{n_{-}}\,\cdot\,P\,\nabla_{X,z}\,u_{\,|_{z=-h_0+B(X)}}\;\;,
        $$
        where $\mathbf{n_{-}}$ denotes the outward normal vector to
        the lower boundary of $\Omega$. For the sake of simplicity, the notation
	$\partial_n$ will always denote the outward conormal derivative associated to
	the elliptic problem under consideration.\\

        \noindent \begin{rmrk}When no confusion can be made,
        we denote $\nabla_X$ by $\nabla$. \end{rmrk}

        \vspace{1em}
	
	\noindent As in \cite{Nicholls,BCL,Lannes} we transform the boundary value
        problem (\ref{General})(\ref{bord}) into a new boundary
        problem defined over the flat band
	$$\mathcal{S} = \{(X,z) \in \xR^{d+1} / X \in \xR^d, -1 < z <
        0\}\;\; . $$ 
        Let S be the following diffeomorphism
	mapping $\mathcal{S}$ to $\Omega$ :
	\begin{equation} \label{S}
	S\;:\;\left(
	\begin{array}{ccl}
	\mathcal{S} & \longrightarrow & \Omega \\
	(X,z) & \longmapsto & s(X,z) = (\eta(X)-B(X)+h_0)\,z+\eta(X)
	\end{array}
	\right)\;\;.
	\end{equation}
	\begin{rmrk}As shown in \cite{Lannes}, a more complex ''regularizing'' diffeomorphism must be
	used instead of $(\mathcal{S})$ to obtain a shard dependence on $\eta$ in terms of regularity, but since
	the trivial diffeomorphism $(\mathcal{S})$ suffices for our present purpose, we use it for the sake of simplicity.
	\end{rmrk}
	
	\vspace{1em}
	
	\noindent Clearly, if $v$ is
        defined over $\Omega$ then $\underline{v} = v \circ
	S$ is defined over $\mathcal{S}$. As a consequence, we can set an equivalent
	problem to (\ref{General})(\ref{bord}) on the flat band
        $\mathcal{S}$ using the following proposition (see \cite{Lannes2}
        for a proof) :
	
	\vspace{1em}
		
	\begin{prpstn} \label{prop1}
	$u$ is solution of (\ref{General})(\ref{bord}) if and only
        if $\underline{u} = u \circ S$ is
	solution of the boundary value problem\\
	\begin{equation} \label{S1}
	-\nabla_{X,z}\,.\,\underline{P}\,\nabla_{X,z}\,\underline{u} = 0 \;\;\;\;\mbox{in}\;\;\; \mathcal{S}\;\;,
	\end{equation}
	\begin{equation} \label{S2}
	\underline{u}_{\,|_{z=0}} = f \;\;\; \mbox{and} \;\;\; \partial_n\,\underline{u}_{\,|_{z=
	-1}} = 0\;\;,
	\end{equation}
        
        \vspace{0.5em}

	\noindent where $\underline{P}(X,z)$ is given by \\	
	$$ \underline{P}(X,z) = \frac{1}{\eta+h_0-B}\;\,M^T\,P\;M\;\;, $$
	\\
	$$\mbox{with} \;\;M(X,z) = 
	\left(
	\begin{array}{cc}
	(\eta+h_0-B) I_{d \times d} & \;-(z+1)\,\nabla \eta 	 	
	+z\,\nabla B \\
	0 & 1
	\end{array}
	\right) .
	$$
	\end{prpstn}

	\vspace{1em}
	
	\noindent Consequently, let us consider boundary value problems
        belonging to the class (\ref{S1})(\ref{S2}). From now
        on, all references to the problem set on $\mathcal{S}$ will be
        labelled with an underscore.\\

        \noindent On the class (\ref{S1})(\ref{S2}) of problems set on the flat band $\mathcal{S}$,
        we have the following classical existence theorem :
		assuming that $\underline{P}$ and all its derivatives are bounded on  $\mathcal{S}$, 
		if $f \in H^{k+\frac{3}{2}}(\xR^d)$ 
		then there exists a unique solution $u \in
        H^{k+2}(\mathcal{S})$ to (\ref{S1})(\ref{S2}).
        The proof is very classical and we omit it here.
        
        \vspace{2em}

        \noindent As previously seen, we need to consider the
        following operator $Z(\eta,B)$ which maps the value of $u$ at
        the upper bound to the value of $\partial_z
        u|_{z=\eta}$ :

        $$
		Z(\eta,B) :
		\left(
		\begin{array}{cll}
		\vspace{0.5em}
		H^{\frac{3}{2}}(\xR^d) & \longrightarrow & H^{\frac{1}{2}}(\xR^d)\\
		\vspace{0.5em}
		f & \longmapsto & \partial_z u|_{z=\eta} \;\;\;\mbox{with $u$ solution of (\ref{General})(\ref{bord})}
		\end{array}
		\right)\;\;.
		$$

		\noindent 
		\begin{rmrk}The operator $Z_{\varepsilon}$ defined in (\ref{Ze}) corresponds to the operator $Z$ in the
		case where $P = \left( \begin{array}{cc} \varepsilon I_{d} & 0 \\ 0 & 1 \end{array} \right)$ in 
		$(\ref{General})(\ref{bord})$.\end{rmrk}
		
        \vspace{1em}

        \noindent To construct an approximation of this operator $Z(\eta,B)$, we
        need the following lemma which gives a coercitivity result taking into account
        the anisotropy of $(\ref{General})(\ref{bord})$.

        \vspace{1em}

	\begin{lmm}\label{lemma1}
	Let $\eta \in W^{1,\infty}(\xR^{d})$ and  $B \in
        W^{1,\infty}(\xR^{d})$. Then for all $V \in \xR^{d+1}\,:$ 
        $$  (V\,,\,\underline{P}\,V) \;\ge\;
        c_0(\,||\eta||_{W^{1,\infty}},||B||_{W^{1,\infty}})\;|\sqrt{P}\,V|^{\,2}\;\;,
	$$
        where $c_0$ is a strictly positive function
        given by
        $$
        c_0(x,y) =
        \frac{h_{min}}{(d+1)^2}\,\min\left(1,\frac{1}{h_{min}(x+h_0+y)},
        \frac{\displaystyle\min_{1 \le i \le d}\frac{p_{d+1}}{p_i}}{
        (x+y)^2}\right)\;\;.
        $$
	\end{lmm}

        \vspace{2em}

        \noindent \begin{proof}
        Using Proposition \ref{prop1} , we can write, with $\delta(X) =
        \eta(X)+h_0-B(X)$ :
        \begin{eqnarray*}
        (V\,,\,\underline{P}\,V)  & = & \Big(\,\frac{1}{\delta}\;V \,,\, M^T\,P\;M\;V\Big)\\
        & = & \Big(\,\frac{1}{\delta}\;M\;V \,,\, P\;M\;V\Big)\\
        & = & \Big(\,\frac{1}{\delta}\;\sqrt{P}\;M\;V \,,\, \sqrt{P}\;M\;V\Big)\\
        & = & \Big|\,\frac{1}{\displaystyle\sqrt{\delta}}\,\mathcal{M}\,(\sqrt{P}\;V)\Big|^{\,2}
        \end{eqnarray*}
        where $\mathcal{M} =
        \sqrt{P}\,M\,(\sqrt{P})^{-1}$. Thanks to the
        condition (\ref{cond}), we deduce the invertibility of $M$
        and hence the invertibility of $\mathcal{M}$. This yields
        the following norm inequality for all $U \in \xR^{d+1}$ :
        $$ |U| \leq (d+1)
        \left|\displaystyle\sqrt{\delta}\;\mathcal{M}^{-1}\right|
        _{\infty}\;\left|\frac{1}{\displaystyle\sqrt{\delta}}\,\mathcal{M}\,U\right| \;\;,$$
        with
        $$\mathcal{M}^{-1} = 
		\left(
		\begin{array}{cc}
			\vspace{0.5em}
		\frac{1}{\delta}\, I_{d \times d} & \;\frac{1}{\delta\sqrt{p_{d+1}}}\,\sqrt{P^d}\,((z+1)\,\nabla \eta 	 	
		-z\,\nabla B) \\
		0 & 1
		\end{array}
		\right) .
		$$
        where $|A|_{\infty} = \displaystyle\sup_{1 \le i,j
          \le d+1} |a_{i,j}|_{L^{\infty}(\xR^{d})}$ and $P^d$ is the $d\times d$ diagonal matrix whose coefficients are
        $(p_i)_{1 \le i \le d}$.\\
        If we apply the previous inequality to our problem, one gets :
        $$ (V\,,\,\underline{P}\,V) \geq
        \displaystyle\frac{1}{(d+1)^2\left|\displaystyle\sqrt{\delta}\;\mathcal{M}^{-1}
        \right|_{\infty}^{\;2}}\;\left|\sqrt{P}\;V\,\right|^{\,2} \;\;.$$
        Thanks to the expression of $\mathcal{M}^{-1}$ given above, we obtain the following inequality
        :
        $$ (V\,,\,\underline{P}\,V) \geq
        c_0(\,||\eta||_{W^{1,\infty}},||B||_{W^{1,\infty}})\;|\sqrt{P}\;V|^2 \;\;,$$
        \vspace{0.5em}where 
        $c_0$ as in the statement of the Lemma \ref{lemma1}. \end{proof}

        \vspace{2em}

        \noindent Let us introduce the space $H^{k,0}(\mathcal{S})$ :
	$$ H^{k,0}(\mathcal{S}) = \{v \in L^2(\mathcal{S}),\;||v||_{H^{k,0}} := \left( \int_{-1}^0 |v(\cdot,z)|_{H^k(\xR^d)}^{\;\;2} dz 	 	
	\right)
	^{\frac{1}{2}} < +\infty \}\;\;. 
        $$

        \vspace{1em}
        
        \noindent The result of this subsection consists in the following theorem which
        aims at giving a rigourous method for deriving an
        asymptotic development of $Z(\eta,B)$. Of course, $P$, and thus $\underline{P}$, 
		as well as the boundaries $\eta$ and $B$, can depend on $\varepsilon$ in the following theorem. In such cases,
		the proof can be easily adapted just by remembering that $0 <
        \varepsilon < 1$.

        \vspace{4em}

        \begin{thrm} \label{prop2}
        Let $p\in\xN^*, k\in\xN^*,\,\eta\in
        W^{k+2,\infty}(\xR^d)$ and $B \in W^{k+2,\infty}(\xR^d)$. Let
        $0 < \varepsilon < 1$ and
        $u_{app}$ be such that
		\begin{equation} \label{Ga1}
		-\nabla_{X,z}\cdot\underline{P}\,\nabla_{X,z}\,u_{app}
                = \varepsilon^p\, R^{\varepsilon} \;\;\;\;\mbox{in}\;\;\;
		\mathcal{S}\;\;,
		\end{equation}	
		\begin{equation} \label{Ga2}
		u_{app\,|_{z=0}} = f \;\;\; , \;\;\;\;
                \partial_n\,u_{app\,|_{z=-1}} =
                \varepsilon^p\,r^{\varepsilon}\;\;,
		\end{equation}
	
			\vspace{0.5em}
	
		\noindent where
                $(R^{\varepsilon})_{\,0<\varepsilon<1}$ and $(r^{\varepsilon})_{\,0<\varepsilon<1}$
                are bounded independently of $\varepsilon$
                respectively in $H^{k+1,0}(\mathcal{S})$ and $H^{k+1}(\xR^d)$.\\
		Assuming that $h_{min}$ is independent of
                $\varepsilon$ and that 
		the coefficients $(p_i)_{1 \le i \le
        d+1}$ of $P$ are such that
        $(\frac{p_i}{p_{d+1}})_{1 \le i \le d}$
        are bounded by a constant $\gamma$ independent of
        $\varepsilon$, we have
        $$
        \left|\,Z(\eta,B)f -\frac{1}{\eta+h_0-B }\,(\partial_z u_{app})_{\,|_{z=0}}\right|_
	{H^{k+\frac{1}{2}}} 
        \le \frac{\varepsilon^p}{\sqrt{p_{d+1}}}\,C_{k+2} \left(||R^{\varepsilon}||_{
	H^{k+1,0}}+|r^{\varepsilon}|_{
	H^{k+1}}\right)\;\;, $$
        where $C_{k+2}=C(|\eta|_{W^{k+2,\infty}},|B|_{W^{k+2,\infty}})$
        and C is a non decreasing function of its arguments,
        independent of the coefficients $(p_i)_{1 \le i \le d+1}$.
	\end{thrm}
	
	\vspace{2em}
	
	\noindent \begin{proof}
        In this proof, we often use the notation $C_k =
        C(|\eta|_{W^{k,\infty}},|B|_{W^{k,\infty}},h_0,h_{min},k,d,\gamma)$
        where C is an undefined non decreasing function of its arguments.
        The notation $C_k$ can thus refer to different constants, but of the same kind.
		\vspace{0.5em} \\
		A simple calculus shows that
        $Z(\eta,B)$
        can be expressed in terms of the
        solution $u$ of (\ref{S1})(\ref{S2}) via the following
        relation :
        $$ Z(\eta,B)f = \frac{1}{\eta+h_0-B} \partial_z \underline{u}_{|_{z=0}} \;\;.$$
	Using this fact, we can write 
	$$ Z(\eta,B)f - \frac{1}{\eta+h_0-B}\;\partial_z u_{app\,|_{z=0}} = 
	\frac{1}{\eta+h_0-B}\;\partial_z (\underline{u} - u_{app})_{\,|_{z=0}} \;\;.$$
	Introducing $\varphi := u_{app} - \underline{u}$ we use a trace
        theorem (see Metivier \cite{Metivier} p.23-27) to get
	\begin{equation} \label{trace}
        |\,Z(\eta,B)f - \frac{1}{\eta+h_0-B}\,
	(\partial_z u_{app})_{\,|_{z=0}} |_
	{H^{k+\frac{1}{2}}} \le C_{k+1} (||\partial_z\varphi||_{H^{k+1,0}}
	+ ||\partial_z^2\varphi||_{H^{k,0}}) \;\;.
	\end{equation}
	It is clear that the proof relies on finding an adequate
        control of $||\partial_z \varphi||_{H^{k+1,0}}$ and
        $||\partial_z^2\varphi||_{H^{k,0}}$.
        The rest of this proof will hence be devoted
        to the estimate of both terms.

        \vspace{1em}

	\noindent {\bf 1. } Let us begin with the estimate of $||\partial_z
        \varphi||_{H^{k+1,0}}$. To deal correctly with this problem, we introduce the following
        norm
        $||.||_{\Dot{H^1}}$ defined by :
	$$ ||\varphi||_{\dot{H^1}} := ||\sqrt{P}\,\nabla_{X,z}\varphi||_{L^2(\mathcal{S})} \;\;.$$
        First remark that for all $\alpha \in \xN^d$ such that
        $|\alpha| \le k,\, \partial^{\alpha}\varphi
	$ solves :
	\begin{equation} \label{dl}
	\left\{
	\begin{array}{l}
	\vspace{0.5em}
	-\nabla_{X,z}\cdot\underline{P}\,\nabla_{X,z}\,\partial^{\alpha}\varphi = \varepsilon^p\,\partial^{\alpha}
	R^{\varepsilon} + \nabla_{X,z}\cdot[\partial^{\alpha},\underline{P}]\,\nabla_{X,z}\,\varphi \;\;,\\
	\partial^{\alpha}\varphi_{\,|z=0} = 0\;, \;\;\;\;\;\;\;
	\partial_n(\partial^{\alpha}\varphi)_{\,|z=-1}
	+ \partial_{\;\;n}^{[\partial^{\alpha},\underline{P}]}\varphi_{\,|z=-1} =
        \varepsilon^p\,\partial^{\alpha}r^{\varepsilon}\;\;.
	\end{array}
	\right.
	\end{equation}
	In order to get an adequate control of the norm $||\partial_z
        \varphi||_{H^{k+1,0}}$, we prove the following
        estimate by induction on $|\alpha| \le k$:
	\begin{equation} \label{major}
	\forall \alpha \in \xN^d\;/\;|\alpha| \le k, \;\;\;\; ||\partial^{\alpha} \varphi||_{\dot{H^1}} 
	\le \frac{\varepsilon^p}{\sqrt{p_{d+1}}}\,C_{k+1} \left(||R^{\varepsilon}||_{
	H^{k,0}}+|r^{\varepsilon}|_{
	H^{k}}\right)\;\;.
	\end{equation}

	\noindent The proof of (\ref{major}) is hence divided into two
        parts : initialization of the induction and heredity.

        \vspace{1em}

        \noindent $\bullet$ Initialization : $|\alpha| = 0\;.$\\
        Taking $\alpha = 0$, multiplying (\ref{dl}) by
        $\varphi$ and integrating by parts leads to :
	$$ \left(\,\underline{P}\,\nabla_{X,z}\,\varphi,\;\nabla_{X,z}\,\varphi\,
	\right)_{L^2(\mathcal{S})} + \int_{\xR^d}\partial_n
        \varphi_{|_{z=0}}\varphi_{|_{z=0}} - \int_{\xR^d}\partial_n
        \varphi_{|_{z=-1}}\varphi_{|_{z=-1}}
	= \left(\,\varepsilon^p\,
	R^{\varepsilon},\;\varphi\,
	\right)_{L^2(\mathcal{S})}\;\;.$$
        The boundary term at the free surface vanish because of the condition
        $\varphi_{\,|z=0} = 0$ and using the condition at the bottom
        leads to :
        $$
        \left(\,\underline{P}\,\nabla_{X,z}\,\varphi,\;\nabla_{X,z}\,\varphi\,
	\right)_{L^2(\mathcal{S})}
	= \left(\,\varepsilon^p\,
	R^{\varepsilon},\;\varphi\,
	\right)_{L^2(\mathcal{S})}+\varepsilon^p \int_{\xR^d}
        r^{\varepsilon}\varphi_{|_{z=-1}}\;\;.
        $$
	Finally, using Cauchy-Schwartz
        inequality, one gets :
	\begin{equation} \label{int}
	\left(\,\underline{P}\,\nabla_{X,z}\,\varphi,\;\nabla_{X,z}\,\varphi\,
	\right)_{L^2(\mathcal{S})} \le \,\varepsilon^p\,||
	R^{\varepsilon}||_{L^2(\mathcal{S})}\,||\varphi||_{L^2(\mathcal{S})}+\varepsilon^p\,|
	r^{\varepsilon}|_{L^2(\xR^d)}\,|\varphi_{|_{z=-1}}|_{L^2(\xR^d)}\;\;.
        \end{equation}
    Recalling that $\varphi_{\,|_{z=0}} = 0$ and that the band $\mathcal{S}$
        is bounded in the vertical direction, one can use Poincar\'e inequality so that
        $||\varphi||_{L^2(\mathcal{S})} \le ||\partial_z \varphi
        ||_{L^2(\mathcal{S})}$ and $|\varphi_{|_{z=-1}}|_{L^2(\xR^d)} \le ||\partial_z \varphi
        ||_{L^2(\mathcal{S})} $. Therefore, 
        (\ref{int}) yields
	\begin{eqnarray} \label{reste}
	\left(\,\underline{P}\,\nabla_{X,z}\,\varphi,\;\nabla_{X,z}\,\varphi\,
	\right)_{L^2(\mathcal{S})}
	& \le & \frac{\varepsilon^p}{\sqrt{p_{d+1}}}\,||
	R^{\varepsilon}||_{L^2(\mathcal{S})}\,||\varphi||_{\dot{H^1}}+\frac{\varepsilon^p}{\sqrt{p_{d+1}}}\,|
	r^{\varepsilon}|_{L^2(\xR^d)}\,||\varphi||_{\dot{H^1}} \;\;.
	\end{eqnarray} 
	Using Lemma \ref{lemma1} to bound $\left(\,\underline{P}\,\nabla_{X,z}\,\varphi,\;\nabla_{X,z}\,\varphi\,
	\right)_{L^2(\mathcal{S})}$ from below, one finally gets :
        $$
         c_0 (\,|\eta|_{W^{1,\infty}},
	|B|_{W^{1,\infty}}) ||
        \varphi||_{\dot{H^1}}^{\;\,2} \le \frac{\varepsilon^p}{\sqrt{p_{d+1}}}\,\left(||
	R^{\varepsilon}||_{H^{0,0}}+|
	r^{\varepsilon}|_{H^{0}}\right)\,||\varphi||_{\dot{H^1}}\;\;.
        $$
        Since $c_0 (\,|\eta|_{W^{1,\infty}},|B|_{W^{1,\infty}})$
        depends only on $h_{min}$, $d$ and $\gamma$ through the quantity
        $\min_{1 \le i \le d} \frac{p_{d+1}}{p_i}$ (by Lemma \ref{lemma1}), and
        since the function $c_0$ is a decreasing function of its
        arguments (again by Lemma \ref{lemma1}),
        we get the following desired estimate :
        $$
        ||\varphi||_{\dot{H^1}} \le
        \frac{\varepsilon^p}{\sqrt{p_{d+1}}}\,
        C_{1} \,\left(||
	R^{\varepsilon}||_{H^{0,0}}+|
	r^{\varepsilon}|_{H^{0}}\right)\;\;,
        $$
        which ends the initialization of the induction.

	\vspace{1em}

        \noindent$\bullet$ Heredity : for $m \in \xN^*$ fixed such that
        $m \le k$, we suppose that (\ref{major}) is verified for all
        $\alpha \in \xN^d$ such that $|\alpha| \le m-1$.\\
        Let $\alpha \in \xN^d$ such that $|\alpha| = m$. Multiplying
        (\ref{dl}) by $\partial^{\alpha} \varphi$ and integrating by
        parts on
        $\mathcal{S}$ leads to :
        $$ \left(\,\underline{P}\,\nabla_{X,z}\,\partial^{\alpha}\varphi,\;\nabla_{X,z}\,\partial^{\alpha}\varphi\,
	\right)_{L^2(\mathcal{S})} + \int_{\xR^d}
        \partial^{\alpha}\varphi_{|z=0}
        \,\partial_n\partial^{\alpha}\varphi_{|z=0} - \int_{\xR^d} \partial^{\alpha}\varphi_{|z=-1} \,\partial_n\partial^{\alpha}\varphi_{|z=-1} 
	= \left(\,\varepsilon^p\,\partial^{\alpha}
	R^{\varepsilon},\;\partial^{\alpha}\varphi\,
	\right)_{L^2(\mathcal{S})} $$ $$ - \left(\,[\partial^{\alpha},\underline{P}]\,\nabla_{X,z}\,\varphi,\;\nabla_{X,z}\,
	\partial^{\alpha}\varphi\,\right)_{L^2(\mathcal{S})} - \int_{\xR^d}\partial^{\alpha}\varphi_{|z=0}
        \,\partial_{\;n}^{[\partial^{\alpha},\underline{P}]}\varphi_{|z=0} +  
        \int_{\xR^d}\partial^{\alpha}\varphi_{|z=-1} \,\partial_{\;n}^{[\partial^{\alpha},\underline{P}]}\varphi_{|z=-1}\;\;.$$
	The boundary terms at ${z=0}$ vanish because
        of the condition $\partial^{\alpha}\varphi_{\,|z=0} = 0$, and
	using the second boundary condition $\partial_n(\partial^{\alpha}\varphi)_{\,|z=-1}
	+ \partial_{\;\;n}^{[\partial^{\alpha},\underline{P}]}\varphi_{\,|z=-1} =
        \varepsilon^p \partial^{\alpha}r^{\varepsilon}$, one gets :
        $$
        \left(\,\underline{P}\,\nabla_{X,z}\,\partial^{\alpha}\varphi,\;\nabla_{X,z}\,\partial^{\alpha}\varphi\,
	\right)_{L^2}
	= \left(\,\varepsilon^p\,\partial^{\alpha}
	R^{\varepsilon},\;\partial^{\alpha}\varphi\,
	\right)_{L^2}  +  
        \,\varepsilon^p\int_{\xR^d}\partial^{\alpha}\varphi_{|z=-1}\,\partial^{\alpha}r^{\varepsilon} 
        - \left(\,[\partial^{\alpha},\underline{P}]\,\nabla_{X,z}\,\varphi,\;\nabla_{X,z}\,
	\partial^{\alpha}\varphi\,\right)_{L^2}\;\;,
        $$
        and with Cauchy-Schwartz :
	\begin{eqnarray*}
	\left(\,\underline{P}\,\nabla_{X,z}\,\partial^{\alpha}\varphi,\;\nabla_{X,z}\,\partial^{\alpha}\varphi\,
	\right)_{L^2(\mathcal{S})} & \le & \,\varepsilon^p\,||\partial^{\alpha}
	R^{\varepsilon}||_{L^2(\mathcal{S})}\,||\partial^{\alpha}\varphi||_{L^2(\mathcal{S})}
        + \varepsilon^p\,|\partial^{\alpha}
	r^{\varepsilon}|_{L^2(\xR^d)}\,|\partial^{\alpha}\varphi_{|_{z=-1}}|_{L^2(\xR^d)}
        \\ & & + \left|\left(\,[\partial^{\alpha},\underline{P}]\,\nabla_{X,z}\,\varphi,\;\nabla_{X,z}\,
	\partial^{\alpha}\varphi\,\right)_{L^2(\mathcal{S})}\right|.
	\end{eqnarray*}
        By using the same method and arguments as in the
        initialization, the following inequality arises :
        \begin{equation} \label{pt}
        c_0 (\,|\eta|_{W^{1,\infty}},
	|B|_{W^{1,\infty}}) ||
        \partial^{\alpha}\varphi||_{\dot{H^1}}^{\;\;2} \le \frac{\varepsilon^p}{\sqrt{p_{d+1}}}\,\left(||
	R^{\varepsilon}||_{H^{k,0}}+|
	r^{\varepsilon}|_{H^{k}}\right)\,||\partial^{\alpha}\varphi||_{\dot{H^1}} + \left|\left(\,[\partial^{\alpha},\underline{P}]\,\nabla_{X,z}\,\varphi,\;\nabla_{X,z}\,
	\partial_{\alpha}\varphi\,\right)_{L^2(\mathcal{S})}\right|.
        \end{equation}
	
	\noindent Let us now focus on the second term of the left hand side of
        (\ref{pt}). In order to get an adequate control of this
        term, we have to write explicitly the commutator
        $[\partial^{\alpha},\underline{P}]$ :
        $$ [\partial^{\alpha},\underline{P}]\,\nabla_{X,z}\,\varphi =
        \sum_{\begin{array}{c}
        \scriptstyle \alpha'+\alpha'' = \alpha\\
        \scriptstyle \alpha' \neq 0
        \end{array}} C(|\alpha'|,|\alpha''|)
        \partial^{\alpha'}\underline{P}\,\nabla_{X,z}\,\partial^{\alpha''}\varphi\;\;,
        $$
        where $C$  is a constant depending only on $|\alpha'|$ and $|\alpha''|$.
        This leads to the expression
        $$
        \left(\,[\partial^{\alpha},\underline{P}]\nabla_{X,z}\,\varphi,\,\nabla_{X,z}\,
	\partial^{\alpha}\varphi\,\right)_{L^2(\mathcal{S})}
        = \sum_{\begin{array}{c}
        \scriptstyle \alpha'+\alpha'' = \alpha\\
        \scriptstyle \alpha' \neq 0
        \end{array}} C(|\alpha'|,|\alpha''|)
        \left(\, \partial^{\alpha'}\underline{P}\,\nabla_{X,z}\,\partial^{\alpha''}\varphi\;,\;\nabla_{X,z}\,
	\partial^{\alpha}\varphi\,\right)_{L^2(\mathcal{S})}\;\;.
        $$
        From now on, we just consider a single term of this sum.
        Using Proposition \ref{prop1} we derive the explicit
        expression of $\underline{P}$ and deduce from it the explicit
        expression of $\partial^{\alpha'}\underline{P}$ :
        $$
        \partial^{\alpha'} \underline{P} = 
        \left(
        \begin{array}{cc}
        \vspace{0.5em}
        \left(\partial^{\alpha'}\eta-\partial^{\alpha'}B\right)\,P_{d} &
        P_{d}\;\partial^{\alpha'}\mathcal{U}\\
        \left(P_{d}\;\partial^{\alpha'}\mathcal{U}\right)^{T} & 
        \partial^{\alpha'}        \left(\frac{p_{d+1}+\mathcal{U} \cdot P_{d}\;\mathcal{U}}{\eta+h_0-B}\right)
        \end{array}
        \right)\;\;,
        $$
        where $P_{d}$ is the diagonal $(d\times d)$ matrix 
        whose coefficents are $(p_i)_{1\le i
          \le d}$, and $\mathcal{U}$ the vector defined by
        $\mathcal{U} = -(z+1)\nabla\eta
        +z\nabla B$. Using this expression, one
        easily gets (with $\nabla = \nabla_X$) :
        \begin{eqnarray}\label{4}
        &\left(\partial^{\alpha'}\underline{P}\,\nabla_{X,z}\,\partial^{\alpha''}\varphi\;,\;\nabla_{X,z}\,
	\partial^{\alpha}\varphi\,\right)_{L^2(\mathcal{S})} =
        \left(\,(\partial^{\alpha'}\eta-\partial^{\alpha'}B)P_{d}\nabla\partial^{\alpha''}\varphi\,,\,\nabla \partial^{\alpha}\varphi\right)_{L^{2}(\mathcal{S})}&\nonumber \\ 
        &+ \left(\,\partial_z \partial^{\alpha''}\varphi
        P_{d}\;\partial^{\alpha'}\mathcal{U}\,,\,\nabla
        \partial^{\alpha} \varphi\right)_{L^{2}(\mathcal{S})}
        + \left(\,P_{d}\;\partial^{\alpha'}\mathcal{U}\cdot \nabla
        \partial^{\alpha''}\varphi\,,\,\partial_z \partial^{\alpha}
        \varphi\right)_{L^{2}(\mathcal{S})}&\nonumber \\ 
        &+ \left(\,\partial^{\alpha'}
        \left(\frac{p_{d+1}+\mathcal{U} \cdot
            P_{d}\;\mathcal{U}}{\eta+h_0-B}\right)\,\partial_z \partial^{\alpha''}
        \varphi\,,\,\partial_z \partial^{\alpha}
        \varphi\right)_{L^{2}(\mathcal{S})}.&
        \end{eqnarray}
        If we focus on the first term of the right hand side of this
        equality, we easily get the following intermediate control
        using Cauchy-Schwartz inequality and the definition of
        $||.||_{\dot{H^{1}}}$ :
        \begin{eqnarray*}
        \left(\,(\partial^{\alpha'}\eta-\partial^{\alpha'}B)P_{d}\,\nabla\partial^{\alpha''}\varphi\,,\,\nabla \partial^{\alpha}\varphi\right)_{L^{2}(\mathcal{S})}
        & \le &
        \left(|\eta|_{W^{|\alpha'|,\infty}}+|B|_{W^{|\alpha'|,\infty}}\right)||\sqrt{P_{d}}\,\nabla\partial^{\alpha''}\varphi||_{L^{2}(\mathcal{S})}
        ||\sqrt{P_{d}}\,\nabla\partial^{\alpha}\varphi||_{L^{2}(\mathcal{S})}\\
        & \le &
        \left(|\eta|_{W^{k,\infty}}+|B|_{W^{k,\infty}}\right)||\partial^{\alpha''}\varphi||_{\dot{H^{1}}}
        ||\partial^{\alpha}\varphi||_{\dot{H^{1}}}\\
        & \le &
        \frac{\varepsilon^p}{\sqrt{p_{d+1}}}\,C_{k+1}\left(||
	R^{\varepsilon}||_{H^{k,0}}+|
	r^{\varepsilon}|_{H^{k}}\right)
        ||\partial^{\alpha}\varphi||_{\dot{H^{1}}}\;\;.
        \end{eqnarray*}
        To derive the last inequality, we used the induction hypothesis on
        $||\partial^{\alpha''}\varphi||_{\dot{H^{1}}}$ since
        $|\alpha''| \le m-1$.  \\
        Let us now focus on the second term of the right hand side of
        (\ref{4}). Using the same arguments as previously and Poincar\'e
        inequality, the following controls arise :
        \begin{eqnarray*}
        \left(\,\partial_z \partial^{\alpha''}\varphi
        P_{d}\;\partial^{\alpha'}\mathcal{U}\,,\,\nabla
        \partial^{\alpha} \varphi\right)_{L^{2}(\mathcal{S})}
        & \le &
        ||\sqrt{P_{d}}\partial^{\alpha'}\mathcal{U}||_{\infty}||\partial_z \partial^{\alpha''}\varphi||_{L^{2}(\mathcal{S})}
        ||\sqrt{P_{d}}\,\nabla\partial^{\alpha}\varphi||_{L^{2}(\mathcal{S})}\\
        & \le &
        \sqrt{\frac{||P_{d}||_{\infty}}{p_{d+1}}}
        \left(|\eta|_{W^{|\alpha'|+1,\infty}}+|B|_{W^{|\alpha'|+1,\infty}}\right)\,
        ||\partial^{\alpha''}\varphi||_{\dot{H^{1}}}
        ||\partial^{\alpha}\varphi||_{\dot{H^{1}}}\\
        & \le &
        \frac{\varepsilon^p}{\sqrt{p_{d+1}}}\,
        C_{k+1}\,
        \left(||
	R^{\varepsilon}||_{H^{k,0}}+|
	r^{\varepsilon}|_{H^{k}}\right)\,
        ||\partial^{\alpha}\varphi||_{\dot{H^{1}}}\;\;,
        \end{eqnarray*}
        since $\frac{||P_{d}||_{\infty}}{p_{d+1}} \le
        \gamma$.\\
        The control of the
        third term of the right hand side of (\ref{4}) comes in the
        same way :
        $$
        \left(\,
        P_{d}\;\partial^{\alpha'}\mathcal{U} \cdot \nabla
        \partial^{\alpha''} \varphi\,,\,\partial_z \partial^{\alpha}\varphi\right)_{L^{2}(\mathcal{S})}
        \le
        \frac{\varepsilon^p}{\sqrt{p_{d+1}}}\,
        C_{k+1}\,
        \left(||
	R^{\varepsilon}||_{H^{k,0}}+|
	r^{\varepsilon}|_{H^{k}}\right)\,
        ||\partial^{\alpha}\varphi||_{\dot{H^{1}}}\;\;.
        $$

        \noindent The next step consists in controling the last
        term of the right hand side of (\ref{4}). We need
        to do some preliminary work on this term before attempting to
        estimate it adequately. A straightforward calculus gives us :
	\begin{eqnarray*}
        \partial^{\alpha'}\left(\frac{p_{d+1}+\mathcal{U}\cdot P_{d}\,\mathcal{U}}{\eta+h_0-B}\right)
        & = & \sum_{\begin{array}{c}
        \scriptstyle \beta_1+\beta_2 = \alpha'\\
        \beta_1 \neq 0
        \end{array}}
        C(|\beta_1|,|\beta_2|)\partial^{\beta_1}\left(\mathcal{U}\cdot P_{d}\,\mathcal{U}\right)
        \partial^{\beta_2}\left(\frac{1}{\eta+h_0-B}\right)\\
        & & +\,
        \mathcal{U}\cdot
        P_{d}\,\mathcal{U}\,\partial^{\alpha'}\left(\frac{1}{\eta+h_0-B}\right) + p_{d+1}\,\partial^{\alpha'}\left(\frac{1}{\eta+h_0-B}\right)\;\;.
        \end{eqnarray*}
        We plug the previous writing and use the
        same tools as previously to get :
        \begin{eqnarray*}
        \left(\,\partial^{\alpha'}
        \left(\frac{p_{d+1}+\mathcal{U} \cdot
        P_{d}\;\mathcal{U}}{\eta+h_0-B}\right)\,\partial_z \partial^{\alpha''}
        \varphi\,,\,\partial_z \partial^{\alpha}
        \varphi\right)_{L^{2}(\mathcal{S})}
        & \le &
        ||P_{d}||_{\infty}\,C_{k+1}\,||\partial_z \partial^{\alpha''}\varphi||_{L^{2}(\mathcal{S})}\,
        ||\partial_z
        \partial^{\alpha}\varphi||_{L^{2}(\mathcal{S})}\\
        & & + \,C_{k+1}\,||\sqrt{p_{d+1}}\,\partial_z \partial^{\alpha''}\varphi||_{L^{2}(\mathcal{S})}\,
        ||\sqrt{p_{d+1}}\,\partial_z
        \partial^{\alpha}\varphi||_{L^{2}(\mathcal{S})}\\
        & \le &
        \frac{||P_{d}||_{\infty}}{p_{d+1}}\,C_{k+1}\,||\partial^{\alpha''}\varphi||_{\dot{H^{1}}}
        ||\partial^{\alpha}\varphi||_{\dot{H^{1}}}\\
        & & + C_{k+1}\,||\partial^{\alpha''}\varphi||_{\dot{H^{1}}}\,
        ||\partial^{\alpha}\varphi||_{\dot{H^{1}}}\\
        & \le &
        \frac{\varepsilon^p}{\sqrt{p_{d+1}}}\,
        C_{k+1}\,\left(||
	R^{\varepsilon}||_{H^{k,0}}+|
	r^{\varepsilon}|_{H^{k}}\right)\,
        ||\partial^{\alpha}\varphi||_{\dot{H^{1}}}\;\;,
        \end{eqnarray*}
        where we once more used the induction hypothesis.\\
        Gathering the four previous estimations of each term of the
        right hand side of (\ref{4}) and using the explicit writing of
        the commutator $[\partial^{\alpha},\underline{P}]$ leads to the final estimate of 
        $\left|([\partial^{\alpha},\underline{P}]\,\nabla_{X,z}\,\varphi\;,\;\nabla_{X,z}\,
	\partial^{\alpha}\varphi\,)_{L^2(\mathcal{S})}\right|$ :
        \begin{flushleft}
        $
        \left|\left([\partial^{\alpha},\underline{P}]\,\nabla_{X,z}\,\varphi\;,\;\nabla_{X,z}\,
	\partial^{\alpha}\varphi\,\right)_{L^2(\mathcal{S})}\right|
        \le
        \frac{\varepsilon^p}{\sqrt{p_{d+1}}}\,
        C_{k+1}\,\left(||
	R^{\varepsilon}||_{H^{k,0}}+|
	r^{\varepsilon}|_{H^{k}}\right)\,
        ||\partial^{\alpha}\varphi||_{\dot{H^{1}}}\;\;.
        $
        \end{flushleft}
        The last step simply consists in pluging this last estimation
        in the estimation (\ref{pt}), which gives :
        \begin{equation} \label{ptf}
        c_0 (\,|\eta|_{W^{1,\infty}},
	|B|_{W^{1,\infty}}) ||
        \partial^{\alpha}\varphi||_{\dot{H^1}}^{\;\;2} \le  \frac{\varepsilon^p}{\sqrt{p_{d+1}}}\,||
	R^{\varepsilon}||_{H^{k,0}}\,||\partial^{\alpha}\varphi||_{\dot{H^1}}\;\;.
        \end{equation}
        As in the initialization, this last estimation leads to the
        desired result, which ends the heredity and hence the proof of
        (\ref{major}).

        \vspace{1em}

        \noindent To conclude this first part of the proof, we use the fact that :
        \begin{eqnarray*}
        ||\partial_z\varphi||_{H^{k+1,0}} & \le & C(k+1)\,\sup_{|\alpha|
          \le k+1}
        ||\partial_z\partial^{\alpha}\varphi||_{L^2(\mathcal{S})} \\
        & \le & \frac{C(k+1)}{\sqrt{p_{d+1}}}\,\sup_{|\alpha|
          \le k+1}
        ||\partial^{\alpha}\varphi||_{\dot{H^1}}\;\;,
        \end{eqnarray*}
        and the estimate (\ref{major}) we just proved to finally get :
        \begin{equation} \label{firstcontrol}
        ||\partial_z\varphi||_{H^{k+1,0}} \le \frac{\varepsilon^p}{p_{d+1}}\,C_{k+2}\,
        \left(||
	R^{\varepsilon}||_{H^{k+1,0}}+|
	r^{\varepsilon}|_{H^{k+1}}\right)\;\;,
        \end{equation}
        which ends the first part of the proof.

        \vspace{2em}

	\noindent {\bf 2. } In this second part, we aim at controlling
        the quantity $||\partial_z^2
        \varphi||_{H^{k,0}}\,$. To this end, we prove with a direct
        method the following estimate :
        \begin{equation} \label{secondcontrol}
        ||\partial_z^2\varphi||_{H^{k,0}} \le
        \frac{\varepsilon^p}{p_{d+1}}\,C_{k+2}\,\left(||
	R^{\varepsilon}||_{H^{k+1,0}}+|
	r^{\varepsilon}|_{H^{k+1}}\right)\;\;.
        \end{equation}
        \noindent We first use the equation satisfied
        by $\varphi$ in order to express $\partial_z^2
        \varphi$ in terms of other derivatives of $\varphi$ such as
        $\nabla \varphi, \partial_z \varphi, \nabla\partial_z\varphi$
        or $\Delta \varphi$. There comes the following expression :
        \begin{eqnarray*}
        \partial_z^2 \varphi &=&  \left(\frac{\eta +
          h_0 -B}{p_{d+1}+ \mathcal{U}\cdot
          P_{d}\,\mathcal{U}}\right)\,[-\nabla_{X,z}\cdot
        Q\,\nabla_{X,z}\varphi
        -\frac{\partial_z\left(\mathcal{U}\cdot P_{d}\,\mathcal{U}\right)}{\eta+h_0-B}\,\partial_z \varphi - \varepsilon^p\,R^{\varepsilon}]\;,
        \end{eqnarray*}
        where 
        $
        Q = \displaystyle\left(
        \begin{array}{cc}
        (\eta+h_0-B)\,P_{d} & P_{d}\,\mathcal{U}\\  
        \left(P_{d}\,\mathcal{U}\right)^T & 0
        \end{array}
        \right)
        $.\\
        \\
        The following estimates arise (using
        $||u||_{H^{k,0}} \le C(k) \sup_{|\alpha| \le k}
        ||\partial^{\alpha} u||_{L^2(\mathcal{S})}\,$) :
        \begin{eqnarray*}
        ||\partial_z^2 \varphi||_{H^{k,0}} & \le &
        \frac{1}{p_{d+1}}\,C_{k}\,\displaystyle\Big[\,||\nabla_{X,z}\cdot
        Q\,\nabla_{X,z}\varphi||_{H^{k,0}}
        +
        C_{k+1}\,||P_{d}||_{\infty}\,||\partial_z
        \varphi||_{H^{k,0}} +
        \varepsilon^p||R^{\varepsilon}||_{H^{k,0}}\;\Big]\;\;, \\
        & \le &
        \frac{1}{p_{d+1}}\,C_{k+1}\,\Big[\;C(k)\,\sup_{|\alpha| \le
            k} ||\partial^{\alpha}\left(\nabla_{X,z}\cdot
        Q\,\nabla_{X,z}\varphi\right)||_{L^2(\mathcal{S})}
        +
        \frac{||P_{d}||_{\infty}}{\sqrt{p_{d+1}}}\,C(k)\,\sup_{|\alpha| \le
          k} ||\partial^{\alpha}\varphi||_{\dot{H^1}}\\
        & & +
        \varepsilon^p||R^{\varepsilon}||_{H^{k,0}}\;\Big]\;\;, \\
        & \le & \frac{1}{p_{d+1}}\,C_{k+1}\;\sup_{|\alpha| \le
          k} ||\partial^{\alpha}\left(\nabla_{X,z}\cdot
        Q\,\nabla_{X,z}\varphi\right)||_{L^2(\mathcal{S})}
        +
        \frac{\varepsilon^p}{p_{d+1}}\,\frac{||P_{d}||_{\infty}}{p_{d+1}}\,C_{k+1}\,\left(||
	R^{\varepsilon}||_{H^{k,0}}+|
	r^{\varepsilon}|_{H^{k}}\right)\\
        & & +
        \frac{\varepsilon^p}{p_{d+1}}\,C_{k+1}\,||R^{\varepsilon}||_{H^{k,0}}\;\;, \\
        & \le & \frac{1}{p_{d+1}}\,C_{k+1}\;\sup_{|\alpha| \le
          k} ||\partial^{\alpha}\left(\nabla_{X,z}\cdot
        Q\,\nabla_{X,z}\varphi\right)||_{L^2(\mathcal{S})}
        +
        \frac{\varepsilon^p}{p_{d+1}}\,C_{k+1}\,\left(||
	R^{\varepsilon}||_{H^{k,0}}+|
	r^{\varepsilon}|_{H^{k}}\right)\;\;,
        \end{eqnarray*}
        where we used the result (\ref{major}) and
        the fact that $\frac{||P_{d}||_{\infty}}{p_{d+1}} \le \gamma$.\\
        The last part of the initialization aims at correctly estimating the norm $||\partial^{\alpha}\left(\nabla_{X,z}\cdot
        Q\,\nabla_{X,z}\varphi\right)||_{L^2(\mathcal{S})}$. The explicit writing
        $$
        \partial^{\alpha}\left(\nabla_{X,z}\cdot
        Q\,\nabla_{X,z}\varphi\right) = \sum_{\alpha'+\alpha''=\alpha} C(|\alpha'|,|\alpha''|)\,\left(\nabla_{X,z}\cdot\partial^{\alpha'}Q\,\nabla_{X,z}\,\partial^{\alpha''}\varphi\right)\;\;,
        $$
        and the expression of $Q$ furnishes us with the following estimates :
        
        \begin{eqnarray} \label{interm}
        ||\nabla_{X,z}\cdot\partial^{\alpha'}Q\,\nabla_{X,z}\,\partial^{\alpha''}\varphi||_{L^2(\mathcal{S})}&\le&
        C_{|\alpha'|}\;||\nabla\cdot
        P_{d} \nabla\partial^{\alpha''}\varphi||_{L^2(\mathcal{S})} \nonumber \\
        & & + C_{|\alpha'|+1}\;
        ||\sqrt{P_{d}}||_{\infty}\,||\sqrt{P_{d}}\,\nabla\partial^{\alpha''}\varphi||_{L^2(\mathcal{S})}\nonumber \\
        & & + C_{|\alpha'|+2}\; ||P_{d}||_{\infty}\,||\partial_z\partial^{\alpha''}\varphi||_{L^2(\mathcal{S})}\nonumber \\
        & & + C_{|\alpha'|+1}\; ||P_{d}||_{\infty}\,||
        \partial_z\nabla\partial^{\alpha''}\varphi||_{L^2(\mathcal{S})}\;)\nonumber \\
        &\le&  C_{k+2}\;(\;||\nabla\cdot
        P_{d} \nabla\partial^{\alpha''}\varphi||_{L^2(\mathcal{S})}\nonumber \\
        & & +
        ||\sqrt{P_{d}}||_{\infty}\,
        ||\partial^{\alpha''}\varphi||_{\dot{H^1}} +
        \frac{||P_{d}||_{\infty}}{\sqrt{p_{d+1}}}
        ||\partial^{\alpha''}\varphi||_{\dot{H^1}}\;)\nonumber \\
        & \le &  C_{k+2}\;\left(\;||\nabla\cdot
        P_{d} \nabla\partial^{\alpha''}\varphi||_{L^2(\mathcal{S})} +
        \,\varepsilon^p\,C_{k+1}\,\left(||
	R^{\varepsilon}||_{H^{k,0}}+|
	r^{\varepsilon}|_{H^{k}}\right)\;\right)\;\;.\nonumber \\
        & & 
        \end{eqnarray}
        We estimate the term $||\nabla\cdot
        P_{d} \nabla\partial^{\alpha''}\varphi||_{L^2(\mathcal{S})}$ using the following technique :
        \begin{eqnarray*}
        ||\nabla\cdot
        P_{d} \nabla\partial^{\alpha''}\varphi||_{L^2(\mathcal{S})} &\le&
        ||\sqrt{P_{d}}||_{\infty}\,\sum_{1 \le i \le d} ||\sqrt{p_i} \partial_{x_i}^2 \partial^{\alpha''}\varphi||_{L^2(\mathcal{S})}\\
        &\le&
        ||\sqrt{P_{d}}||_{\infty}\,\sum_{1 \le i \le d} ||\partial_{x_i}\partial^{\alpha''}\varphi||_{\dot{H^1}}\\
        &\le& d\,||\sqrt{P_{d}}||_{\infty}\,\sup_{|m|=|\alpha''|+1}
        ||\partial^{m}\varphi||_{\dot{H^1}} \\
        & \le & d\,||\sqrt{P_{d}}||_{\infty}\,\sup_{|m| \le k+1}
        ||\partial^{m}\varphi||_{\dot{H^1}}\\
        & \le & \varepsilon^p\,C_{k+2}\,\left(||
	R^{\varepsilon}||_{H^{k+1,0}}+|
	r^{\varepsilon}|_{H^{k+1}}\right)\;\;.
        \end{eqnarray*}
        We plug this result in (\ref{interm}) to obtain 
        $$
        ||\nabla_{X,z}\cdot
        \partial^{\alpha'}Q\,\nabla_{X,z}\partial^{\alpha''}\varphi||_{L^2(\mathcal{S})} \le
        \varepsilon^p\,C_{k+2}\,\left(||
	R^{\varepsilon}||_{H^{k+1,0}}+|
	r^{\varepsilon}|_{H^{k+1}}\right)\;\;,
        $$
        which finally leads to
        $$
        ||\partial^{\alpha}\left(\nabla_{X,z}\cdot
        Q\,\nabla_{X,z}\varphi\right)||_{L^2(\mathcal{S})} \le
        \varepsilon^p\,C_{k+2}\,\left(||
	R^{\varepsilon}||_{H^{k+1,0}}+|
	r^{\varepsilon}|_{H^{k+1}}\right)\;\;.
        $$
        This way, we get our desired estimation of $||\partial_z^2 \varphi||_{H^{k,0}}$: 
        \begin{equation*}
        ||\partial_z^2 \varphi||_{H^{k,0}} \le 
        \frac{\varepsilon^p}{p_{d+1}}\,C_{k+2}\,\left(||
	R^{\varepsilon}||_{H^{k+1,0}}+|
	r^{\varepsilon}|_{H^{k+1}}\right)\;\;.
        \end{equation*}

        \vspace{1em}

        \noindent Gathering (\ref{firstcontrol}) and (\ref{secondcontrol}) in (\ref{trace})
        ends the proof of the theorem.
        \end{proof}

        \vspace{1em}

	\subsection{Application}

	\noindent We recall that by definition, 
        $Z_{\varepsilon}(\varepsilon\eta,\beta b)f = \partial_z u_{|_{z=\varepsilon
            \eta}}$ where $u$ is solution of the boundary value problem
	\begin{equation} \label{Ge1}
	\varepsilon\,\Delta u + \partial_z^2 u = 0 \;\;\; \mbox{in}\; \Omega\;\;,
	\end{equation}	
	\begin{equation} \label{Ge2}
	u_{|_{z=\varepsilon \eta}} = f  \;\;\;,\;\;\;
        \left(\partial_z u - \varepsilon \beta \nabla b \cdot \nabla u\right)_{|_{z=-1+\beta
	b}} = 0\;\;,
	\end{equation}	
        
        \vspace{1em}

	\noindent This elliptic problem
        (\ref{Ge1})(\ref{Ge2}) belongs to the class of
        general elliptic problems (\ref{General})(\ref{bord}) defined
        in the previous subsection. The
        corresponding matrix $P$ is here designed by
        $P^{\varepsilon}$ :
        \\
		\begin{equation} \label{Pe}
		P^{\varepsilon} = \left(
		\begin{array}{cc}
		\varepsilon\,I_{d \times d} & 0 \\
		0 & 1
		\end{array}
		\right)\;\;.
		\end{equation}

        \vspace{1em}

	\noindent The upper boundary of $\Omega$ is here defined by $\{z = \varepsilon \eta\}$ and
        the lower one by 
	$\{z = -1 + \beta b\}$. We make the additionnal assumption that
        $\varepsilon$ and $\beta$ are bounded in the following sense : $0 <
        \varepsilon < 1$ and
        there exists a strictly positive constant $\beta_0$ such that  and $0 < \beta < \beta_0$.
        Furthermore, condition
	(\ref{cond}) is here verified thanks to condition
        (\ref{hmin}). And finally, we remark that
        $(\frac{p_i}{p_{d+1}})_{1\le i \le d}$ are bounded by 1 since $0 <
        \varepsilon < 1$. Our goal is here to apply the previous theorem 
	to get asymptotic estimates on $Z_{\varepsilon}(\varepsilon\eta,\beta b)$. \\

        \vspace{0.5em}

        \noindent We recall that we are here interested in
        two differerent regimes depending on the $\beta$ parameter. 
        The first one, namely $\beta = O(\varepsilon)$, refers to
        the physical case of a bottom with variations of slow
        amplitude.
        The second one, namely $\beta = O(1)$, refers on the contrary
        to variations of high amplitude of the bottom. In order to
        improve the readability, we take $\beta_0 = 1$ : we thus write $\beta = \varepsilon$ for the first regime
        and $\beta = 1$ for the second one. 

        \vspace{1em}        

        \subsubsection{The regime $\beta = \varepsilon$ :
        small variations of bottom topography}

        The boundaries
        of the domain $\Omega$ are here defined by $\{z = \varepsilon \eta\}$ and $\{z =
        -1 + \varepsilon b\}$ while matrix $P^{\varepsilon}$ remains as in (\ref{Pe}).
		Thanks to Proposition \ref{prop1} we are able
        to set an equivalent problem to (\ref{Ge1})(\ref{Ge2}) defined
        over the flat band $\mathcal{S}$ : $\underline{u} = u \circ S$ then solves
        the problem :
        
        \begin{equation} \label{Se1}
		-\nabla_{X,z}\,\cdot\,\underline{P}^{\varepsilon}\,\nabla_{X,z}\,\underline{u} = 0 \;\;\; \mbox{in}\;\; \mathcal{S}\;\;,
		\end{equation}	
		\begin{equation} \label{Se2}
		\underline{u}_{|_{z=0}} = f  \;\;\;,\;\;\; \partial_n\,\underline{u}_{|_{z=-1}} = 0\;\;.
		\end{equation}	
        where the matrix $\underline{P}^{\varepsilon}$ is given by
        $$
        \underline{P}^{\varepsilon} = \left(
		\begin{array}{cc}
			\vspace{1em}
		\varepsilon(1+\varepsilon(\eta- b))\,I_{d \times d} &
			-\varepsilon^2[(z+1)\nabla \eta -  z \nabla b] \\
		 -\varepsilon^2[(z+1)\nabla \eta -  z \nabla b]^T & \frac{1+\varepsilon^3\left|(z+1)\nabla \eta -  
		z \nabla b\right|^2}{1+\varepsilon(\eta- b)}
		\end{array}
		\right)\;\;.
        $$
        The
        following result gives a rigourously justified asymptotic
        expansion
        of $Z_{\varepsilon}(\varepsilon\eta,\beta b)f$ as $\varepsilon$ goes to $0$ :
        
        \begin{prpstn} \label{prop21}
        \vspace{1em}
        Let $k \in \xN,\, \eta \in W^{k+2,\infty}(\xR^d)$ and $b \in
        W^{k+2,\infty}(\xR^d)$.\\ Then for all $f$ such that $\nabla f
        \in H^{k+6}(\xR^d)$, we have :
        $$
        \left| Z_{\varepsilon}(\varepsilon\eta,\beta b)f - (\varepsilon Z_1 +
          \varepsilon^2 Z_2) \right|_{H^{k+1/2}} \le
          \varepsilon^3
          C_{k+2}\,|\nabla f|_{H^{k+6}}\,,
        $$
        with :
        $$
        \left\{
        \begin{array}{l}
        \vspace{0.5em} Z_1 := - \Delta f\;\;, \\ Z_2 :=
        -\frac{1}{3}\Delta^2\,f - (\eta- b)\Delta f  +
        \nabla b \cdot \nabla f\;.
        \end{array}
        \right.
        $$
        \end{prpstn}

        \vspace{2em}

        \noindent \begin{proof}
        To prove this proposition, we use essentially Theorem
        \ref{prop2} with $p=3$. We know that
        $\left(\frac{p_i}{p_{d+1}}\right)_{1 \le i \le d}$ are bounded
        by $1$. Thus, in order to derive an asymptotic
        expansion of $Z_{\varepsilon}(\varepsilon\eta,\beta b)f$, we only need to
        compute an approximate solution $u_{app}$ which satisfies
        the hypothesis of Theorem \ref{prop2} for $p=3$. This approximate
        solution $u_{app}$ can be constructed as in \cite{BCL} using a classical WKB method,
        which consists in looking for $u_{app}$ under the form $u_{app}=
        u_0 + \varepsilon u_1 + \varepsilon^2 u_2$. We want this
        function to verify the properties required by Theorem \ref{prop2},
        that is to say :
        \begin{equation} \label{Ga11}
		-\nabla_{X,z}\,.\,\underline{P}^{\varepsilon}\,\nabla_{X,z}\,u_{app} = \varepsilon^p\,R^{\varepsilon} 
		\;\;\;\;\mbox{in}\;\;\;
		\mathcal{S}\;\;,
		\end{equation}
		\begin{equation} \label{Ga12}
		u_{app\,|_{z=0}} = f \;\;\; , \;\;\;\;
                \partial_n\,u_{app\,|_{z=-1}} = \varepsilon^p r^{\varepsilon}\;\;.
		\end{equation}
        
        \vspace{0.5em}        
        
        \noindent where $(R^{\varepsilon})_{\,0<\varepsilon<1}$ and $(r^{\varepsilon})_{\,0<\varepsilon<1}$
        are bounded independently of $\varepsilon$ respectively in $H^{k+1,0}(\mathcal{S})$ and $H^{k+1}(\xR^d)$.\\
        We decompose the matrix $\underline{P^{\varepsilon}}$ under the
        form $\underline{P^{\varepsilon}} = P_0 + \varepsilon P_1 +
        \varepsilon^2 P_2 + \varepsilon^3 P_{\varepsilon}$ where $P_0,
        P_1, P_2$ are independent of $\varepsilon$, and
        if we plug the desired expression of $u_{app}$ into this
        problem,
        we get $R^{\varepsilon} = \nabla\cdot
        T^{\varepsilon}$ and $r^{\varepsilon} = \mathbf{e_z} \cdot
        T^{\varepsilon}_{|_{z=-1}}$ where  $T^{\varepsilon} = P_2 \nabla_{X,z} u_1 + P_1
        \nabla_{X,z} u_2 + P_{\varepsilon} \nabla_{X,z} (u_o+u_1+u_2)$,
        and the following system of equations and boundary
        conditions on
        $u_0,u_1,u_2$ :
        $$
        \left\{
        \begin{array}{l}
        \vspace{0.5em}
        \partial_z^2 u_0 = 0\;\;,\\
        \vspace{0.5em}
        \partial_z^2 u_1 + \left(\Delta - (\eta -  b)
          \partial_z^2 \right) u_0 = 0\;\;,\\
          \vspace{0.5em}
        \partial_z^2 u_2 + \left(\Delta - (\eta -  b)
          \partial_z^2 \right) u_1
        +(\eta- b) \Delta u_0 - 2\left[(z+1)\nabla f -
           z \nabla b \right]\cdot \\\vspace{0.5em} \nabla \partial_z u_0 
		- \left[(z+1)\Delta f -
          z \Delta b \right]\cdot \partial_z u_0 - (\eta -
         b)^2 \partial_z^2 u_0 = 0\;\;,
        \end{array}
        \right.
        $$
        $$
        \;\;\mbox{with}\;\;\left\{
        \begin{array}{l}
        \vspace{0.5em}
        u_{0|_{z=0}} = f\;\;,\\
        \vspace{0.5em}
        u_{i|_{z=0}} = 0 ,\;\;\; 1 \le i \le 2\;\;,\\
        \vspace{0.5em}
        \partial_z u_{i|_{z=-1}} = 0 ,\;\;\; 0 \le i \le 1\;\;,\\
        \partial_z u_{2|_{z=-1}} -  \nabla b \cdot \nabla u_{0|_{z=-1}} = 0\;\;.
        \end{array}
        \right.
        $$
        
        \vspace{1em}

        \noindent We can verify that the following values of
        $u_0,u_1,u_2$ satisfy the previous equations and boundary
        conditions :
        \begin{eqnarray*}
        u_0 & = & f\;\;,\\
        u_1 & = & \left(\frac{1}{2}-\frac{(z+1)^2}{2}\right)\Delta f\;\;,\\
        u_2 & = & \left(\frac{(z+1)^4}{24} - \frac{(z+1)^2}{4} +
          \frac{5}{24} \right) \Delta^2 f + \left(1-(z+1)^2\right) (\eta -
         b) \Delta f +  z \nabla b \cdot \nabla f\;\;.
        \end{eqnarray*}
        Using these values of $u_0,u_1$ and $u_2$, one can easily
        obtain the following estimates of
        $R^{\varepsilon}$ and $r^{\varepsilon}$ :
        $$
        ||R^{\varepsilon}||_{H^{k+1,0}} \le C_{k+2}\,|\nabla
        f|_{H^{k+6}}\;\;,\;\;\; |r^{\varepsilon}|_{H^{k+1}} \le C_{k+2}\,|\nabla
        f|_{H^{k+3}}\;\;.
        $$
        Thus $u_{app}$
        satisfies the properties required to apply Theorem
        \ref{prop2}. The last steps of the proof consists in
        computing $\frac{1}{1+\varepsilon(\eta-
          b)}\partial_z u_{app|_{z=0}}$ using the explicit expression
        of $u_{app}$ previously determined, and then apply
        Theorem \ref{prop2}.
        An easy Taylor expansion of $\frac{1}{1+\varepsilon(\eta-
          b)}\partial_z u_{app|_{z=0}}$ then yields the result.
        \end{proof}

        \vspace{2em}

        \begin{rmrk}The method developped here to get and
        prove our asymptotic expansion is improved compared to the
        one developped in BCL since we do not need
        here to compute the term $u_3$.\end{rmrk}
        
        \begin{rmrk}
        If we take $b=0$ - i.e. if we consider a flat bottom - , we of
        course get the same
        expansion as the ones proved in \cite{BCL}.
        \end{rmrk}

        \vspace{1em}
        
	\subsubsection{The regime $\beta = 1$ : strong variations
        of bottom topography}

        The boundaries of the domain $\Omega$
        are here defined by $\{z = \varepsilon \eta\}$ and $\{z =
        -1 +  b\}$. Using again Proposition \ref{prop1},
        we 
        set an equivalent problem to (\ref{Ge1})(\ref{Ge2}) defined
        over the flat band $\mathcal{S}$ : this new problem is the
        same
        as the one defined in the first regime, at the exception of
        the matrix $\underline{P}^{\varepsilon}$ which is now given by
        $$
        \underline{P}^{\varepsilon} = \left(
		\begin{array}{cc}
			\vspace{1em}
		\varepsilon(1+\varepsilon \eta -  b)\,I_{d \times d} &
			-\varepsilon [\varepsilon (z+1) \nabla \eta -  z \nabla b] \\
		 -\varepsilon [\varepsilon (z+1) \nabla \eta -  z
			 \nabla b]^T & \frac{1+\varepsilon\left|\varepsilon (z+1)
				 \nabla \eta -  z \nabla b \right|^2}{1+\varepsilon
			   \eta- b}
		\end{array}
		\right) \;\;.
        $$
        As in the first regime we give a rigourously justified asymptotic expansion of
        $Z_{\varepsilon}(\varepsilon\eta,\beta b)f$
        in the present regime.

        \begin{prpstn} \label{prop22}
        \vspace{1em}
        Let $k \in \xN,\, \eta \in W^{k+2,\infty}(\xR^d)$ and $b \in
        W^{k+2,\infty}(\xR^d)$.\\ Then for all $f$ such that $\nabla f
        \in H^{k+6}(\xR^d)$, we have :
        $$
        \left| Z_{\varepsilon}(\varepsilon\eta,\beta b)f - (\varepsilon Z_1 +
          \varepsilon^2 Z_2) \right|_{H^{k+1/2}} \le
          \varepsilon^3
          C_{k+2}\,|\nabla f|_{H^{k+6}}\,,
        $$
        with :
        $$
        \left\{
        \begin{array}{l}
        \vspace{0.5em} Z_1 := - \nabla \cdot \Big((1-b)\,\nabla f\Big) \;\;,\\ Z_2 :=
        \frac{1}{2}\,\nabla \cdot \left(\frac{1}{3} (1-b)^3\,\nabla \Delta f -
        (1-b)^2\,\nabla \nabla \cdot \Big((1-b)\,\nabla f\Big)\right) - \eta \Delta f \;.
        \end{array}
        \right.
        $$
        \end{prpstn}

        \vspace{0.5em}

        \noindent \begin{proof}
        The proof of this proposition follows exactly the same steps as
        the proof of Proposition \ref{prop21}. The following values of $u_0,
        u_1, u_2$ are found :
        \begin{eqnarray*}
        u_0 & = & f\;\;,\\
        u_1 & = & \frac{(1-b)^2}{2}(1-(z+1)^2)\Delta f + z\,(1-b)\nabla
        b\cdot\nabla f\;\;,\\
        u_2 & = & \frac{(1-b)^4}{24}\,\Delta^2 f\,z^4 +
        \frac{(1-b)^3}{6}\,\Delta\nabla\cdot\Big((1-b)\nabla f\Big)\,z^3 - (1-b) \eta\Delta
        f\,z^2 \\ & & +
        \Big[\frac{(1-b)}{2}\,\nabla\cdot\left(\frac{(1-b)^3}{3}\nabla\Delta f -
           (1-b)^2 \nabla\nabla\cdot\Big((1-b)\nabla f\Big)\right) \\  & & -\eta\Big(2(1-b)\Delta f
        +\nabla b \cdot \nabla f\Big)\Big]\,z\;\;.
        \end{eqnarray*}
        The error bound is derived in the same way and the previous
        values lead to the result.
        \end{proof}

        \vspace{1em}

        \begin{rmrk} By formally taking $b=\varepsilon b$, we recover
          the result of Proposition \ref{prop21}.
        \end{rmrk}
	
        \vspace{1em}

        \subsection{Derivation of Boussinesq-like models for uneven bottoms}

        \noindent We recall the Zakharov formulation of the
        water waves equations, from which we intend to derive new systems
        using results of the previous subsection :
	$$
		(S_{0})\left\{
		\begin{array}{l}
		\vspace{1em}
		\partial_t \psi - \varepsilon \partial_t \eta Z_{\varepsilon}(\varepsilon\eta,\beta b)\psi
		+\frac{1}{2}\,\Big[\,\varepsilon
		\,|\nabla \psi - \varepsilon\,\nabla \eta Z_{\varepsilon}(\varepsilon\eta,\beta b)\psi |^2
		+ \,|Z_{\varepsilon}(\varepsilon\eta,\beta b)\psi |^2
		\, \Big] + \eta = 0 \;\;,\\
		\partial_t \eta + \varepsilon \nabla \eta \cdot \Big[\,\nabla \psi - \varepsilon \nabla 
		\eta Z_{\varepsilon}(\varepsilon\eta,\beta b)
		\psi \,\Big] = \frac{1}{\varepsilon}\,Z_{\varepsilon}(\varepsilon\eta,\beta b)\psi\;\;.
		\end{array}
		\right.
        $$
        As in \cite{BCL}, we introduce the
        notion of consistency.
        \begin{dfntn}
        Let $\sigma, s \in \xR$, $\varepsilon_0 > 0$, $T>0$ and let ${(V^{\varepsilon},
        \eta^{\varepsilon})_{0<\varepsilon<\varepsilon_0}}$ be
		bounded in
        $W^{1,\infty}([0,\frac{T}{\varepsilon}];\,H^{\sigma}(\xR^d)^{d+1})$ independently of $\varepsilon$. 
        This family
        is called consistent (at order $s$) with a system $(S)$ if
        it is solution of $(S)$ with a residual of order $\varepsilon^2$ in $L^{\infty} 
        ([0,\frac{T}{\varepsilon}];\,H^{s}(\xR^d)^{d+1})$. 
        \end{dfntn}
        
        \noindent Thanks to the previous results,
        we are now able to enounce the following
        propositions which show the consistency of two
        Boussinesq-like systems with the system $(S_{0})$. We introduce here the quantity
        $h = 1-b$ which corresponds to the non-dimensional still water depth. From now on,
        this quantity is considered as a bottom term since it only depends on $b$.
        
        \vspace{1em}

        \begin{prpstn}[Small variations regime $\beta = \varepsilon$] \label{Bouss1}
		Let $T>0$, $s \ge 0$, $\sigma \ge s$ and $(\psi^{\varepsilon},
          \eta^{\varepsilon})_{0<\varepsilon<\varepsilon_0}$ be a family
          of solutions
          of (\ref{S0}) such that $(\nabla\psi^{\varepsilon},
          \eta^{\varepsilon})_{0<\varepsilon<\varepsilon_0}$ is bounded with respect to $\varepsilon$ in
          $W^{1,\infty}([0,\frac{T}{\varepsilon}];\,H^{\sigma}(\xR^d)^{d+1})$ with $\sigma$ large enough. 
          We define $V^{\varepsilon} :=
          \nabla\psi^{\varepsilon}$ . Then the family ${(V^{\varepsilon},
          \eta)_{0<\varepsilon<\varepsilon_0}}$ is consistent with
          the following system :
          $$
          (\mathcal{B}_1)\left\{
          \begin{array}{l}
          \vspace{1em} \partial_t V + \nabla \eta +
          \frac{\varepsilon}{2} \nabla |V|^2 = 0 \;\;,\\
          \partial_t \eta + \nabla \cdot V+
          \varepsilon\left[\nabla\cdot\Big((\eta-b)\,V\Big) +
          \frac{1}{3}\,\Delta\nabla\cdot V\right] = 0\;\;.
          \end{array}
          \right.
          $$ 
       \end{prpstn}

       \noindent \begin{proof}
       This is clear thanks to the asymptotic expansion
       of the operator $Z_{\varepsilon}(\varepsilon\eta,\beta b)$ : plugging this in
       system (\ref{S0}), neglecting the terms of order
       $O(\varepsilon^2)$, and taking the gradient yields the result. \end{proof}

       \vspace{1em}

       \begin{prpstn}[Strong variations regime $\beta = 1$] \label{Bouss2}
         Let $T>0$, $s \ge 0$, $\sigma \ge s$ and $(\psi^{\varepsilon},
          \eta^{\varepsilon})_{0<\varepsilon<\varepsilon_0}$ be a family
          of solutions
          of (\ref{S0}) such that $(\nabla\psi^{\varepsilon},
          \eta^{\varepsilon})_{0<\varepsilon<\varepsilon_0}$ is bounded with respect to $\varepsilon$ in
          $W^{1,\infty}([0,\frac{T}{\varepsilon}];\,H^{\sigma}(\xR^d)^{d+1})$ with $\sigma$ large enough. 
          We define $V^{\varepsilon} :=
          \nabla\psi^{\varepsilon}$. Then the family ${(V^{\varepsilon},
          \eta)_{0<\varepsilon<\varepsilon_0}}$ is consistent with
          the following system (with $h = 1- b$) :
          $$
          (\mathcal{B}_2)\left\{
          \begin{array}{l}
          \vspace{1em} \partial_t V + \nabla\eta +
          \frac{\varepsilon}{2} \nabla |V|^2 = 0\;\;,\\
          \partial_t \eta + \nabla \cdot (h\,V) +
          \varepsilon\left[\nabla\cdot (\eta V) -
          \frac{1}{2}\nabla\cdot\left(\frac{h^3}{3}\nabla\nabla\cdot
          V - h^2 \nabla \nabla \cdot (h\,V)\right)\right] = 0\;\;.
          \end{array}
          \right.
          $$ 
       \end{prpstn}
        
       \vspace{1em} 
        
       \noindent These results close the first part. From now on,
       our work is divided into two parts, each corresponding to one of
       our two regimes $\beta = \varepsilon$ and $\beta = 1 $.
	   Each analysis starts from the two previous
       Boussinesq-like models. 
	
        \vspace{2em}
	
	\section{The regime of small
         topography variations}
	
        \vspace{1em}

       We recall our previously derived Boussinesq-like system $(\mathcal{B}_{1})$ on which we base our
       analysis :
       $$
       (\mathcal{B}_1)\left\{
         \begin{array}{l}
           \vspace{1em} \partial_t V + \nabla \eta +
           \frac{\varepsilon}{2} \nabla |V|^2 = 0 \;\;,\\
           \partial_t \eta + \nabla \cdot V+
           \varepsilon\left[\nabla\cdot\Big((\eta-b)\,V\Big) +
             \frac{1}{3}\,\Delta\nabla\cdot V\right] = 0\;\;.
         \end{array}
       \right.
       $$

       \noindent We now follow the method
       put forward in \cite{BCS} and \cite{BCL} to derive equivalent systems to $(\mathcal{B}_1)$
       in the meaning of consistency. The rigourous justifications of the derivation of these systems is 
       adressed in section 2 in \cite{BCL}. 

       \subsection{A first class of equivalent systems}

       As in \cite{BCS} (for the 1D case) and in \cite{BCL}, we define :
       $$
       V_{\theta} = \left(1+\frac{\varepsilon}{2}(1-\theta^2)\Delta\right)\,V\;\;,
       $$
       which corresponds to the approximation at the order $\varepsilon^2$ of
       the horizontal component of the velocity field at height
       $-1+\theta $ for $ \theta \in [0,1]$. 
       If we remark that $V_{\theta} =
       \left(1+\frac{\varepsilon}{2}(\theta^2-1)\Delta\right)^{-1}\,V + O(\varepsilon^2)$, the expression of $V_{\theta}$ in terms of
       $V$ comes in the following way by supposing $V$ regular enough :
       $$
       V =
       \left(1+\frac{\varepsilon}{2}(\theta^2-1)\Delta\right)\,V_{\theta} + O(\varepsilon^2)\;\;,
       $$
       where $O(\varepsilon^2)$ is to be taken in the
       $L^{\infty}([0,\frac{T}{\varepsilon}], H^s(\xR^d))$
       norm. \\
       Plugging this relation into the system $(\mathcal{B}_1)$
       leads to :
       $$
       \left\{
         \begin{array}{l}
           \partial_t V_{\theta} + \nabla \eta +
           \frac{\varepsilon}{2} \Big( \nabla |V|^2 +
             (\theta^2-1)\Delta \partial_t V_{\theta} \Big) = O(\varepsilon^2)\;\;,\\
           \partial_t \eta + \nabla \cdot V_{\theta}+
           \varepsilon\left[\nabla\cdot\Big((\eta-b)\,V\Big) +
             \left(\frac{\theta^2}{2}-\frac{1}{6}\right)\,\Delta\nabla\cdot V_{\theta}\right] = O(\varepsilon)^2\;\;.
         \end{array}
       \right.
       $$
       At this point we use the classical BBM trick which consists in writing
       the following approximations at order $O(1)$ , coming from
       the previous equations :
       $$
       \begin{array}{l}
       \vspace{1em}
       \partial_t V_{\theta} = - \nabla \eta +
       O(\varepsilon) = (1-\mu) \partial_t V_{\theta} - \mu
       \nabla \eta + O(\varepsilon)\;\;,\\
       \nabla \cdot V_{\theta} = - \partial_t \eta + O(\varepsilon) = \lambda \nabla \cdot V_{\theta} - (1-\lambda)
       \partial_t \eta + O(\varepsilon)\;\;,
       \end{array}
       $$
       where $\lambda$ and $\mu$ are two real parameters.\\
       We plug these relations into the dispersive terms of the last system to
       get :
       $$
       \left\{
       \begin{array}{l}
       \vspace{1em} \partial_t V_{\theta} + \nabla \eta + \frac{\varepsilon}{2} \left[
         \nabla|V_{\theta}|^2 - \mu (\theta^2-1) \Delta \nabla \eta - (\mu-1)(\theta^2-1) \Delta\partial_t
         V_{\theta} \right] = O(\varepsilon^2)\;\;,\\
       \partial_t \eta + \nabla \cdot V_{\theta} + \varepsilon \left[
         \nabla \cdot\Big((\eta-b) V_{\theta}\Big) + \lambda \left(\frac{\theta^2}{2} -
         \frac{1}{6}\right)\Delta\nabla \cdot V_{\theta} - (1-\lambda)\left(\frac{\theta^2}{2} -
         \frac{1}{6}\right)\Delta\partial_t \eta
       \right] = O(\varepsilon^2)\;\;.
       \end{array}
       \right.
       $$

       \vspace{1em}
       
       \noindent We then introduce the class $\mathcal{S}$ of all the systems of the
       previous form : these systems are denoted by
       $S_{\theta,\lambda,\mu}^{\;\;1}$ and can be rewritten in the compact
       form :
       $$
       (S_{\theta,\lambda,\mu}^{\;\;\;1})\left\{
         \begin{array}{l}
           \vspace{1em}
           (1-\varepsilon a_2 \Delta)\partial_t V_{\theta} + \nabla \eta +
           \varepsilon\,\left[ \frac{1}{2} \nabla |V_{\theta}|^2 + a_1
             \Delta\nabla \eta \right]
           = 0 \;\;,\\
           (1-\varepsilon a_4\Delta)\partial_t  \eta + \nabla \cdot V_{\theta}+
           \varepsilon\Big[\nabla\cdot\Big( (\eta-b)\,V_{\theta}\Big) + a_3
             \Delta\nabla\cdot V_{\theta}
           \Big] = 0\;\;.
         \end{array}
       \right.
       $$ 
       with
       $$
       \begin{array}{c}
       \vspace{0.5em}
       a_1 = -\mu\frac{\theta^2-1}{2} , \;\;\;a_2 =
         (\mu-1)\frac{\theta^2-1}{2} \;\;,\\
       a_3 = \lambda \left(\frac{\theta^2}{2} -
         \frac{1}{6}\right) , \;\;\;a_4 = (1-\lambda)\left(\frac{\theta^2}{2} -
         \frac{1}{6}\right)\;\;.
       \end{array}
       $$

       \vspace{1em}
       
       \noindent On this class $\mathcal{S}$, the previous computations
       give us the following two results of
       consistency :

       \vspace{0.5em}

       \begin{prpstn}\label{vtheta1}
         Let $\theta \in [0,1]$  and $(\psi^{\varepsilon},
         \eta^{\varepsilon})_{0<\varepsilon<\varepsilon_0}$ be a family of
         solutions of (\ref{S0}) such that $(\nabla\psi^{\varepsilon},
          \eta^{\varepsilon})_{0<\varepsilon<\varepsilon_0}$ is bounded with respect to $\varepsilon$ in
        $W^{1,\infty}([0,\frac{T}{\varepsilon}];\,H^{\sigma}(\xR^d)^{d+1})$ with $\sigma$ large enough.
        We define $V^{\varepsilon} = \nabla
         \psi^{\varepsilon}$ and $V^{\varepsilon}_{\theta} =
         \left(1+\frac{\varepsilon}{2}(1-\theta^2)\Delta\right)\,V^{\varepsilon}$.
         Then for all $(\lambda, \mu) \in \xR^2$, the family $(V^{\varepsilon}_{\theta},
         \eta^{\varepsilon})_{0<\varepsilon<\varepsilon_0}$ is consistent
         with the system $(S_{\theta,\lambda,\mu}^{\;\;\;1})$.
       \end{prpstn}

       \noindent \begin{proof}
       We saw in the previous section that if $(\psi^{\varepsilon},
         \eta^{\varepsilon})_{0<\varepsilon<\varepsilon_0}$ is a
         family of solutions of (\ref{S0}), then
          the family $(\nabla\psi^{\varepsilon},
          \eta^{\varepsilon})_{0<\varepsilon<\varepsilon_0}$ is
         consistent with the system $(\mathcal{B}_1)$. Thanks to the
         previous computations, and since the choice of the parameters $(\lambda, \mu)$ is totally
       free, it is clear that $(V_{\theta}^{\varepsilon},
         \eta^{\varepsilon})_{0<\varepsilon<\varepsilon_0}$ is
         consistent with any system
         $(S_{\theta,\lambda,\mu}^{\;\;\;1})$. \end{proof}

       \vspace{1em}

       \begin{prpstn}
         Up to a change of variables, all the systems belonging to
         the class $\mathcal{S}$ are equivalent in the meaning of consistency.
       \end{prpstn}

       \noindent \begin{proof}
       Let $(\theta, \lambda, \mu) \in [0,1]\times\xR^2$ and
       $(V_{\theta}^{\varepsilon}, \eta^{\varepsilon})_{0<\varepsilon < \varepsilon_0}$ a family
       consistent with $(S_{\theta,\lambda,\mu}^{\;\;\;1})$. 
       We then define, for $\theta_1 \in [0,1]$,
       $$
       V_{\theta_1}^{\varepsilon} = \left(1+\frac{\varepsilon}{2}(1-\theta_{1}^2)\Delta\right)\,\left(1-\frac{\varepsilon}{2}(1-\theta^2)\Delta\right)\,V_{\theta}^{\varepsilon}\;\;;
       $$
       using the fact that
       $\left(1-\frac{\varepsilon}{2}(1-\theta^2)\Delta\right)\,V_{\theta}^{\varepsilon} = V^{\varepsilon} +
       O(\varepsilon^2)$ and the previous proposition,  we easily deduce that the family
       ${(V_{\theta_{1}}^{\varepsilon}, \eta^{\varepsilon})}_{0<\varepsilon < \varepsilon_0}$ is
       consistent with any system $(S_{\theta_{1},\lambda_1,\mu_1}^{\;\;\;1})$
       for any $(\lambda_1,\mu_1) \in \xR^2$. \end{proof}

       \vspace{3em}

       \begin{rmrks}
		\hspace{3em}
       \begin{itemize}
       \item[$\bullet$] By taking $\theta=1, \lambda=1, \mu=0$, we
         remark that the previously derived Boussinesq-like system
         $\mathcal{B}_1$ is actually a member of the class $\mathcal{S}$
       \item[$\bullet$] By taking $\lambda = \mu = 1/2$ and $\theta^2
         = 2/3$, we get $a_1 = a_2 = a_3 = a_4 = \frac{1}{12}$, so that
         the dispersive part of the correponding system
         $(S_{\theta,\lambda,\mu}^{\;\;\;1})$ is
         symmetric. However, the nonlinear terms, that are not affected by
         the choice of $\theta, \lambda, \mu$, are not symmetric :
         this problem is adressed in the next section.
       \item[$\bullet$] In \cite{Chen}, Chen formally studied in 1D the case of
         slowly variating bottoms and derived the same
         class of systems at the exception that she
         considered time-dependent bottoms : her systems thus contain
         additionnal time derivative terms on the bottom that does not
         appear here but could be easily obtained for a time dependent
         bottom.
       \end{itemize}
       \end{rmrks}

       \vspace{1em}

       \subsection{A second class of equivalent systems}
       
       Adapting the nonlinear change of variables of \cite{BCL} to the
       present case of varying depth, we introduce $\tilde{V}$ :
       $$ 
       \tilde{V} = \left(1+\frac{\varepsilon}{2}\,(\eta-b)\right)\,V\;\;.
       $$
       This nonlinear change of variable symetrizes the nonlinear part of the equations.
       \vspace{0.5em}

       \noindent This change of variables only affects the nonlinear terms and
       not the dispersive terms. If
       $(V^{\varepsilon}, \eta^{\varepsilon})_{0<\varepsilon < \varepsilon_0}$
       is consistent with a system $(S_{\theta,\lambda,\mu}^{\;\;\;1})$
       of the class $\mathcal{S}$, then $\tilde{V}^{\varepsilon} =
       (1+\frac{\varepsilon}{2}\,\eta^{\varepsilon})\,V^{\varepsilon} $
       and $\eta^{\varepsilon}$ satisfy the following equations :
       $$
       \left\{
       \begin{array}{l}
       \vspace{1em} 
       (1-\varepsilon a_2 \Delta)\,\partial_t \tilde{V}^{\varepsilon} + \nabla \eta^{\varepsilon} + 
		\varepsilon \Big[\frac{1}{4}
         \nabla |\eta^{\varepsilon}|^{2} + \frac{1}{2}\nabla {|\tilde{V}^{\varepsilon}|}^2 +
         \frac{1}{2}\tilde{V}^{\varepsilon}\,\nabla\cdot \tilde{V}^{\varepsilon} 
         -\frac{1}{2} b\,\nabla \eta^{\varepsilon}+ a_1\Delta\nabla \eta^{\varepsilon} \Big]
       = O(\varepsilon^2)\;, \\
       (1-\varepsilon a_4 \Delta)\,\partial_t \eta^{\varepsilon} + \nabla \cdot \tilde{V}^{\varepsilon} + \varepsilon \left[
         \frac{1}{2}\nabla\cdot\left((\eta^{\varepsilon}-b)\,\tilde{V}^{\varepsilon}\right) + a_3 \Delta \nabla \cdot 
		\tilde{V}^{\varepsilon}  \right] = O(\varepsilon^2)\;\;.
       \end{array}
       \right.
       $$
       As observed in \cite{BCL}, if we consider a two-dimensional domain, that is to say $d=1$,
       the nonlinear terms are actually symmetric. But this is
       not the case in a three-dimensional domain. However we
       can deal with this problem for $d=2$ using the following remark coming from \cite{Lannes2}
       :
       $$
       \frac{1}{2}\nabla {|\tilde{V}^{\varepsilon}|}^2 =
       \frac{1}{4}\nabla {|\tilde{V}^{\varepsilon}|}^2 +
       \frac{1}{2}(\tilde{V}^{\varepsilon}\cdot\nabla)\tilde{V}^{\varepsilon} + \frac{1}{2} \tilde{V}^{\varepsilon} \wedge (\nabla\times \tilde{V}^{\varepsilon})
       $$
       Assuming that $\nabla \times \tilde{V}^{\varepsilon} = O(\varepsilon)$,
       one formally derives the following system :
       $$
       \left\{
       \begin{array}{l}
       \vspace{0.5em} 
       (1-\varepsilon a_2 \Delta)\,\partial_t \tilde{V}^{\varepsilon} + \nabla \eta^{\varepsilon} + 
		\varepsilon \left[\frac{1}{4}
         \nabla |\eta^{\varepsilon}|^{2} + \frac{1}{4}\nabla
         {|\tilde{V}^{\varepsilon}|}^2 
         +\frac{1}{2}(\tilde{V}^{\varepsilon}\cdot\nabla)\tilde{V}^{\varepsilon}
          \right. \\
        \vspace{0.5em} \hspace{12em}\;
        +\frac{1}{2}\,\tilde{V}^{\varepsilon}\,\nabla\cdot
        \tilde{V}^{\varepsilon} -\frac{1}{2}\,b\,\nabla \eta^{\varepsilon}
       + a_1\Delta\nabla \eta^{\varepsilon} \Big]
       = O(\varepsilon^2)\;, \\
       (1-\varepsilon a_4 \Delta)\,\partial_t \eta^{\varepsilon} + \nabla \cdot \tilde{V}^{\varepsilon} + \varepsilon \left[
         \frac{1}{2}\nabla\cdot\left((\eta^{\varepsilon}-b)\,\tilde{V}^{\varepsilon}\right) + a_3 \Delta \nabla \cdot 
		\tilde{V}^{\varepsilon}  \right] = O(\varepsilon^2)\;\;.
       \end{array}
       \right.
       $$
       The nonlinear terms of the previous system are now symmetric
       regardless of the dimension.
       This previous computations are
       summed up in the following
       proposition :
       
       \vspace{1em}

       \begin{prpstn}\label{T1}
         Let $(V^{\varepsilon},
         \eta^{\varepsilon})_{0<\varepsilon<\varepsilon_0}$ be a family
         consistent with a system $(S_{\theta,\lambda,\mu}^{\;\;\;1})$ and
         $\tilde{V^{\varepsilon}} =
         \left(1+\frac{\varepsilon}{2}\,\eta^{\varepsilon}\right)\,V^{\varepsilon}$. 
         If $\nabla \times \tilde{V}^{\varepsilon} =
         O(\varepsilon)$, then the family $(\tilde{V^{\varepsilon}},
         \eta^{\varepsilon})_{0<\varepsilon<\varepsilon_0}$ is consistent with
         the following system :
         $$
         (T_{\theta,\lambda,\mu}^{\;\;\;1})\left\{
             \begin{array}{l}
               \vspace{1em} 
               (1-\varepsilon a_2 \Delta)\,\partial_t V + \nabla \eta + 
		\varepsilon \Big[\frac{1}{4}
         \nabla |\eta|^{2} + \frac{1}{4}\nabla
         {|V|}^2 
         +\frac{1}{2}(V\cdot\nabla)V
          +\frac{1}{2}\,V\,\nabla\cdot
        V -\frac{1}{2}\,b\,\nabla \eta
       + a_1\Delta\nabla \eta \Big]
       = 0\;, \\
       (1-\varepsilon a_4 \Delta)\,\partial_t \eta + \nabla \cdot V + \varepsilon \left[
         \frac{1}{2}\nabla\cdot\Big((\eta-b)\,V\Big) + a_3 \Delta \nabla \cdot 
		V  \right] = 0\;\;.
       \end{array}
       \right.
       $$
       \end{prpstn}

       \vspace{1em}   

       \noindent We introduce the class $\mathcal{T}$ composed with the systems of
       the form $(T_{\theta,\lambda,\mu}^{\;\;\;1})$ for any
       $(\theta,\lambda,\mu) \in [0,1]\times\xR^2$. Using this result, we prove the following proposition :

       \vspace{1em} 

       \begin{prpstn}\label{T2}
         Let $\theta \in [0,1]$  and $(\psi^{\varepsilon},
         \eta^{\varepsilon})_{0<\varepsilon<\varepsilon_0}$ be a family of
         solutions of (\ref{S0}) such that $(\nabla\psi^{\varepsilon},
          \eta^{\varepsilon})_{0<\varepsilon<\varepsilon_0}$ is bounded with respect to $\varepsilon$ in
        $W^{1,\infty}([0,\frac{T}{\varepsilon}];\,H^{\sigma}(\xR^d)^{d+1})$ with $\sigma$ large enough.\\
        We define $\tilde{V}^{\varepsilon} =
         \left(1+\frac{\varepsilon}{2}(\eta-b) \right)\left(1+\frac{\varepsilon}{2}(1-\theta^2)\Delta\right)\,\nabla
         \psi^{\varepsilon}$.
         Then for all $(\lambda, \mu) \in \xR^2$, the family $(\tilde{V}^{\varepsilon},
         \eta^{\varepsilon})_{0<\varepsilon<\varepsilon_0}$ is consistent
         with the system $(T_{\theta,\lambda,\mu}^{\;\;\;1})$.
       \end{prpstn}
       \begin{proof}
       Thanks to Proposition \ref{vtheta1}, the family $(V_{\theta}^{\varepsilon},
         \eta^{\varepsilon})_{0<\varepsilon<\varepsilon_0}$ is
         consistent with the system
         $(S_{\theta,\lambda,\mu}^{\;\;\;1})$ for all $(\lambda, \mu)
         \in \xR^2$, where $V_{\theta}^{\varepsilon} =
         \left(1+\frac{\varepsilon}{2}(1-\theta^2)\Delta\right)\,V^{\varepsilon}$ and $V^{\varepsilon} = \nabla
         \psi^{\varepsilon}$.
         We then use the following
         remark : by hypothesis, the velocity field $V$ is
         irrotationnal, thus $V_{\theta}^{\varepsilon} =
         O(\varepsilon)$ and $\tilde{V}^{\varepsilon} =
         O(\varepsilon)$. Applying the previous proposition yields the result.        
       \end{proof}

       \vspace{1em}

       \subsection{A new class of completely symmetric systems}

       We remarked in the first section that there exists values of
       $(\theta,\lambda,\mu)$ such that the dispersive terms are
       symmetric. Consequently, the corresponding system $(T_{\theta,\lambda,\mu}^{\;\;\;1})$ 
       of the class $\mathcal{T}$ is completely symmetric since both its
       dispersive terms and nonlinear terms are symmetric. We thus introduce
	   the non-empty subclass of $\mathcal{T}$ denoted by $\Sigma$ composed with the systems of the form
       $(T_{\theta,\lambda,\mu}^{\;\;\;1})$ for which we have $a_1 =
       a_3$, $a_2 \ge 0, a_4 \ge 0$. The first condition $a_1 =
       a_3$ symetrizes the nonlinear terms and the last ones $a_2 \ge 0, a_4 \ge 0$ ensure the well-posedness
       of these completely symmetric systems.
	   Indeed, one of the great advantages of these
       systems belonging to $\Sigma$ is that we have a well-posedness over a long time scale :

       \begin{prpstn}\label{symsmall}
       Let $s > \frac{d}{2}+1$ and $(\theta, \lambda, \mu)$ be such that
         the system ($T_{\theta,\lambda,\mu}^{\;\;\;1}$) belongs to the
         class $\Sigma$.\\
       Then for all $(V_0, \eta_0) \in H^s(\xR^d)^{d+1}$, there exists a
       time $T_0 \ge 0$ independent of $\varepsilon$ and a unique solution
       $(V,\eta) \in C([0,\frac{T_0}{\varepsilon}]; H^s(\xR^d)^{d+1})
       \cap \;C^1([0,\frac{T_0}{\varepsilon}]; 
	   H^{s-3}(\xR^d)^{d+1}) $ to
       the system ($T_{\theta,\lambda,\mu}^{\;\;\;1}$) such that
       $(V,\eta)_{|_{t=0}} = (V_0, \eta_0)$. \\
       Furthermore, this unique solution is bounded independently of
       $\varepsilon$ in the following sense : there exists a constant
       $C_0$ independent of $\varepsilon$ such that for all $k$
       verifying $s-3k > \frac{d}{2}+1$, we have :
       $$
       |(V,\eta)|_{W^{k,\infty}([0,\frac{T_0}{\varepsilon}];H^{s-3k}(\xR^d)^{d+1})} \le C_0\;\;.
       $$
       \end{prpstn}

       \vspace{1em}

       \noindent \begin{proof}
       This theorem is a very classical result on hyperbolic symmetric
       quasilinear systems, and we omit the proof here. \end{proof}

       \vspace{2em}

       \noindent As in \cite{BCL}, 
       we are now able to rigorously construct approximate solutions
       to the water waves problem from the solutions of any of these
       symmetric systems. \\
       More precisely, let us consider a solution $(\psi^{\varepsilon},
       \eta^{\varepsilon})$ to the initial system (\ref{S0}) with
       initial data $(\psi_0^{\varepsilon}, \eta_0^{\varepsilon})$
       such that $(\nabla\psi_0^{\varepsilon}, \eta_0^{\varepsilon}) \in
       H^s(\xR^d)^{d+1}$ for a suitably large value of $s$. We define
       $V^{\varepsilon} = \nabla \psi^{\varepsilon}$ and
       $V^{\varepsilon}_0 = \nabla \psi^{\varepsilon}_0$. From this
       solution of the water waves problem, we construct an approximate
       solution as follows :

       \vspace{1.5em}

       \begin{itemize}
       \item[$\bullet$] We first construct what we call here
         approximate initial data, by applying the two successive
         changes of variable on the data $(V_0^{\varepsilon}, \eta_0^{\varepsilon})$ :
       $$
       \left\{
       \begin{array}{l}
         \vspace{0.5em}
       V_{\Sigma,0}^{\varepsilon} = \left(1+\frac{\varepsilon}{2}(
         \eta^{\varepsilon}_0-b)\right)\,\left(
         1+\frac{\varepsilon}{2}(1-\theta^2)\Delta\, \right) \, V^{\varepsilon}_0\;\;,\\
       \eta_{\Sigma,0}^{\varepsilon} = \eta^{\varepsilon}_0\;\;. 
       \end{array}
       \right.
       $$
       \item[$\bullet$] We then choose the parameters $(\theta, \lambda, \mu) \in
         [0,1]\times \xR^2$ such that the system
         $(T_{\theta,\lambda,\mu}^{\;\;\;1})$ belongs to the class
         $\Sigma$ of completely symmetric systems (this choice is
         always possible as we saw previously). Using Proposition \ref{symsmall},
         we know that there exists a unique solution to this system
         with initial data $(V_{\Sigma,0}^{\varepsilon},
         \eta_{\Sigma,0}^{\varepsilon})$ : we denote this solution by
         $(V_{\Sigma}^{\varepsilon}, \eta_{\Sigma}^{\varepsilon})$.
       \item[$\bullet$] From this exact solution of the symmetric system
         $(T_{\theta,\lambda,\mu}^{\;\;\;1})$, we finally construct an approximate solution of the water
         waves problem by successively and approximatively inverting the two changes of
         variable as shown below :
         $$ \label{approx}
       \left\{
       \begin{array}{rcl}
         \vspace{0.5em}
       V_{app}^{\varepsilon} &=& \left(
         1-\frac{\varepsilon}{2}(1-\theta^2)\Delta\, \right) \, \left[
         \Big(1-\frac{\varepsilon}{2}(\eta_{\Sigma}^{\varepsilon}-b)\Big)\,V_{\Sigma}^{\varepsilon} \right]\;\;,\\
       \eta_{app}^{\varepsilon} &=& \eta_{\Sigma}^{\varepsilon}\;\;.
       \end{array}
       \right.
       $$
       \end{itemize}

       \vspace{0.5em}

       \noindent This formal construction of an approximate solution founds its
       mathematical justification in the following theorem which is
       the last result of this section.

       \vspace{2em}

       \begin{thrm}\label{small}
         Let $T_1 \ge 0$, $s \ge \frac{d}{2}+1$, $\sigma \ge s+3$ and
         $(\nabla\psi_0^{\varepsilon},\eta_0^{\varepsilon})$ be in $H^{\sigma}(\xR^d)^{d+1}$.
         Let $(\psi^{\varepsilon},
         \eta^{\varepsilon})_{0 < \varepsilon < \varepsilon_0}$ be a
         family of
         solutions of (\ref{S0}) with initial data
         $(\nabla\psi_0^{\varepsilon},\eta_0^{\varepsilon})$ and such that $(\nabla\psi^{\varepsilon},
          \eta^{\varepsilon})_{0<\varepsilon<\varepsilon_0}$ is bounded with
         respect to $\varepsilon$ in
         $W^{1,\infty}([0,\frac{T_1}{\varepsilon}]; H^{\sigma}(\xR^d)^{d+1})$.
         We define $V^{\varepsilon} = \nabla\psi^{\varepsilon}$ and choose
         $(\theta, \lambda, \mu) \in [0,1]\times\xR^2$  such that
         the system $(T_{\theta,\lambda,\mu}^{\;\;\;1}) \in \Sigma$.
         Then for all $\varepsilon < \varepsilon_0$, there exists $T
         \le T_1$ such that we have :
         $$
         \forall t \in [0,\frac{T}{\varepsilon}]\;,\;\;\;|V^{\varepsilon} -
         V^{\varepsilon}_{app}|_{L^{\infty}([0,t];H^s)} +
         |\eta^{\varepsilon} -
         \eta^{\varepsilon}_{app}|_{L^{\infty}([0,t];H^s)} \le
         C\,\varepsilon^2 t\;\;.
         $$
       \end{thrm}
       
       \vspace{1em}

       \noindent \begin{proof}
       We follow in this proof the strategy put forward in \cite{BCL}:
       estimates are done on
       the symmetric system that provides the approximate
       solution rather than on the initial system (\ref{S0}).\\
       To this end, we take  $(\theta, \lambda, \mu) \in [0,1]\times\xR^2$  such that
         the system $T_{\theta,\lambda,\mu}^{\;\;\;1}$ belongs to the
         class $\Sigma$ of completely symmetric systems.\\
       Since $(\psi^{\varepsilon},
         \eta^{\varepsilon})_{0 < \varepsilon < \varepsilon_0}$ is a
         family of solutions of (\ref{S0}) such that $(\nabla\psi^{\varepsilon},
          \eta^{\varepsilon})_{0<\varepsilon<\varepsilon_0}$ is bounded with
         respect to $\varepsilon$ in
         $W^{1,\infty}([0,\frac{T_1}{\varepsilon}]; H^{\sigma}(\xR^d)^{d+1})$,  using Proposition \ref{vtheta1} implies that $(V_{\theta}^{\varepsilon},
         \eta^{\varepsilon})_{0 < \varepsilon < \varepsilon_0}$ 
         where 
         $
         V^{\varepsilon}_{\theta} =
         \left(1+\frac{\varepsilon}{2}(1-\theta^2)\Delta\right)\,V^{\varepsilon}$
         is consistent with the system
         $S_{\theta,\lambda,\mu}^{\;\;\;1}$.\\
         Moreover, Proposition \ref{T1} states that any family $(V^{\varepsilon},
         \eta^{\varepsilon})_{0 < \varepsilon < \varepsilon_0}$
         consistent with the system $S_{\theta,\lambda,\mu}^{\;\;\;1}$
         is, up to the aforementionned nonlinear change of variables, consistent with
         the system $T_{\theta,\lambda,\mu}^{\;\;\;1}$. Applying this
         result to $(V_{\theta}^{\varepsilon},
         \eta^{\varepsilon})_{0 < \varepsilon < \varepsilon_0}$  shows
		 that the family $(\tilde{V}_{\theta}^{\varepsilon},
         \eta^{\varepsilon})_{0 < \varepsilon < \varepsilon_0}$, where
         $
         \tilde{V}^{\varepsilon}_{\theta} =
         \left(1+\frac{\varepsilon}{2} (\eta^{\varepsilon}-b)\right)\,V_{\theta}^{\varepsilon}\;\;,$
         is actually consistent with the symmetric system
         $T_{\theta,\lambda,\mu}^{\;\;\;1}$.\\
       Thanks to Proposition \ref{symsmall}, we know that there exists a time
       $T_0$ such that there exists a unique
       solution $(V_{\Sigma}^{\varepsilon}, \eta_{\Sigma}^{\varepsilon})$
       to this system with initial data $(V_{\Sigma,0}^{\varepsilon},
       \eta_{\Sigma,0}^{\varepsilon})$ (defined in the previous formal
       construction of the approximate solution). 
       We are now interested in computing the error estimates between $(\tilde{V}_{\theta}^{\varepsilon},
         \eta^{\varepsilon})$
       and $(V_{\Sigma}^{\varepsilon},
       \eta_{\Sigma}^{\varepsilon})$. To this end we define $V =
       \tilde{V}_{\theta}^{\varepsilon} - V_{\Sigma}^{\varepsilon}$
       and $\eta = \eta^{\varepsilon} - \eta_{\Sigma}^{\varepsilon}$.
       Writing the equations satisfied by $V$ and $\eta$ and performing
       standard energy estimates on it leads to the following estimate
       :
       $$
       \forall t \in [0,\frac{T_1}{\varepsilon}]\;,\;\;\;|V|_{L^{\infty}([0,t];H^s)} +
         |\eta|_{L^{\infty}([0,t];H^s)} \le
         C\,\varepsilon^2 t
       $$
       where $T = \min(T_0,T_1)$.
       Inverting the nonlinear change of variables and the
       pseudo-differential one yields the final result. \end{proof}

       \vspace{2em}

       \begin{rmrks}
       \hspace{7em}
       \begin{itemize}
       \item[$\bullet$] The construction of the approximated
         solution of the water waves problem relies on the choice of
         the three parameters $\theta, \lambda, \mu$ such that the
         system $T_{\theta,\lambda,\mu}^{\;\;\;1}$ is completely
         symmetric. A great advantage of this method is that this
         choice is totally free : we are indeed allowed to choose any
         suitable triplet $(\theta, \lambda, \mu)$ we want, and
         contruct our approximate solution from the exact solution of
         the system $T_{\theta,\lambda,\mu}^{\;\;\;1}$. In other
         words, approximate solutions of the water waves problem can be
         constructed starting from the exact solution of {\bf any
           symmetric system} of the class $\Sigma$.
       \item[$\bullet$] Our theorem relies implicitly on the existence of a family $(\psi^{\varepsilon},
         \eta^{\varepsilon})_{0 < \varepsilon < \varepsilon_0}$ of
         solutions to the water waves problem in Sobolev spaces : in 1D-surface, this
         existence have been already proved by Craig \cite{Craig}, Schneider-Wayne \cite{SW} 
		thus this implicit hypothesis of existence of solutions
         is actually a fact. However, in 2D-surface, we have no
         existence result for the water waves problem on a long time
         scale. Lannes proved reccenlty in \cite{Lannes} the existence of solutions to
         this problem in Sobolev spaces in 2D-surface, but we do not
         know if these solutions persist on a long time
         scale. Consequently, the analysis is not totally complete in 2D-surface. 
       \end{itemize}
       \end{rmrks}

       \vspace{2em}

       \section{The regime of strong
         topography variations}
	
        \vspace{1em}

       In this section, attention is given to the regime of strong
       variations of the bottom topography : these variations are here
       of order $O(1)$. We recall the Boussinesq-like system
       $(\mathcal{B}_2)$ derived in the previous section, and the
       fact that solutions of the water waves problem are consistent
       with this system.

       $$
       (\mathcal{B}_2)\left\{
         \begin{array}{l}
           \vspace{1em}\partial_t V + \nabla\eta +
           \frac{\varepsilon}{2} \nabla |V|^2 = 0\;\;,\\\vspace{0.5em}
           \partial_t \eta + \nabla \cdot (h\,V) +
           \varepsilon\Big[\nabla\cdot (\eta V) -
             \frac{1}{2}\nabla\cdot\left(\frac{h^3}{3}\nabla\nabla\cdot
               V - h^2 \nabla \nabla \cdot (h\,V)\right)\Big]= 0\;\;.
         \end{array}
       \right.
       $$ 
       
       \noindent Like in the previous section, we aim at deriving asymptotic
       models, constructing approximate solutions of the water waves
       problem, and justifying these approximations. However, the method
       introduced in the previous section cannot be applyied in the
       exact same way : this regime is indeed much more complex since the bottom
       terms have here a greater influence than in the first
       regime. These bottom terms introduce new difficulties which
       compell us to revise and adapt our strategy.  

       \subsection{A first equivalent system}

       First remark that the bottom term $h$ (recall that $h=1-b$ is the non-dimensional still water depth)
	   appears in the
       first order term of the second equation of $(\mathcal{B}_2)$ whereas 
	   it is not present in the first one
	   : we have here a dissymetry of the order one terms. This
       fact becomes important when it comes to the BBM trick which is unlikely
	   to symetrize these terms. To correctly deal with this regime, we
       have to invert the order of the change of variables, and proceed
       with an adapted nonlinear change of variables first that symetrizes both
       order one terms and non-linear terms.

       \vspace{1em}

       \noindent Taking into account the fact that we have to symmetrize both
       order one terms and nonlinear terms, we introduce the
       following change of variables :

       $$
       \tilde{V} = \left(\sqrt{h} +
         \frac{\varepsilon}{2}\,\frac{\eta}{\sqrt{h}}\right)\, V\;\;.
       $$

       \noindent so that

       $$
       V =  \left(\frac{1}{\sqrt{h}} -
         \frac{\varepsilon}{2}\,\frac{\eta}{h\sqrt{h}}\right)\,
       \tilde{V} + O(\varepsilon^2)\;\;.
       $$

       \noindent Assuming that $\nabla \times \tilde{V}
       = O(\varepsilon)$, we formally
       derive the following system of equations satisfied by $\tilde{V}$ and $\eta$ :

       $$
       \left\{
         \begin{array}{c}
           \vspace{0.5em}
           \partial_t \tilde{V} +\sqrt{h} \nabla \eta  + \frac{\varepsilon}{2\sqrt{h}}
           \left[\frac{1}{4} \nabla
             \eta^2 + \frac{1}{2}\nabla |\tilde{V}|^2 +  (\tilde{V} \cdot \nabla) \tilde{V} +
             \tilde{V} \nabla \cdot \tilde{V}\right.\\
           \vspace{1em}
           \left.+
           \frac{1}{h}\left(\frac{1}{2}(\nabla
             h \cdot \tilde{V}) \tilde{V} -
           |\tilde{V}|^2 \nabla h\right)\right] = O(\varepsilon^2)\;\;,\\
           \vspace{0.5em}
           \partial_t \eta + \nabla (\sqrt{h}\cdot \tilde{V}) +
           \frac{\varepsilon}{2\sqrt{h}} \left[\nabla \cdot (\eta \tilde{V}) -\sqrt{h}\nabla\cdot\left(\frac{h^3}{3}\nabla\nabla\cdot
               (\frac{\tilde{V}}{\sqrt{h}})\right.\right.\\ \left.\left. - h^2 \nabla \nabla \cdot (\sqrt{h}\,\tilde{V})\right) -
           \frac{\eta}{2h}\nabla h \cdot \tilde{V}\right] = O(\varepsilon^2)\;\;.
         \end{array}
       \right.
       $$

       \vspace{1em}

       \noindent We introduce the system $(\mathcal{T}_b)$ that corresponds to the
       homogeneous version of the previous system :

	   $$
       (\mathcal{T}_b)\left\{
         \begin{array}{l}
           \vspace{0.5em}
           \partial_t V +\sqrt{h}\;\nabla \eta  + \frac{\varepsilon}{2}\,
           F_h\left(\begin{array}{cc} V \\ \eta \end{array}\right) = 0\;\;,\\
           \vspace{0.5em}
           \partial_t \eta + \nabla (\sqrt{h}\cdot V) +
           \frac{\varepsilon}{2} \Big[\,f_h\left(\begin{array}{cc} V \\ \eta \end{array}\right) 
		-\nabla\cdot\Big(\frac{h^3}{3}\nabla\nabla\cdot
               (\frac{V}{\sqrt{h}}) - h^2 \nabla \nabla \cdot (\sqrt{h}\,V)\Big)\Big] = 0\;\;,
         \end{array}
       \right.
       $$
	where 
	$$
	\left\{
	\begin{array}{lcl}   
	F_h\left(\begin{array}{cc} V \\ \eta \end{array}\right) &=& \frac{1}{\sqrt{h}}\Big(
	 \frac{1}{4} \nabla \eta^2 + \frac{1}{2}\nabla |V|^2 +  (V \cdot \nabla) V +
     V \nabla \cdot V  + \frac{1}{h}\left(\frac{1}{2}(\nabla h \cdot V) V -
    |V|^2 \nabla h\right)\Big)\;\;, \\
	f_h\left(\begin{array}{cc} V \\ \eta \end{array}\right) &=& \frac{1}{\sqrt{h}}\Big(
	\nabla\cdot(\eta V) - \frac{\eta}{2h} \nabla h \cdot V\Big)\;\;.
	\end{array}
	\right.
	$$
	On this new system $(\mathcal{T}_b)$, we have the following result of consistency :

    \begin{prpstn}
    Consider a family $(\psi^{\varepsilon}, \eta^{\varepsilon})_{0
      < \varepsilon < \varepsilon_0}$ of solutions of (\ref{S0}) such that $(\nabla\psi^{\varepsilon},
          \eta^{\varepsilon})_{0<\varepsilon<\varepsilon_0}$ is bounded with respect to $\varepsilon$ in
        $W^{1,\infty}([0,\frac{T}{\varepsilon}];\,H^{\sigma}(\xR^d)^{d+1})$ with $\sigma$ large enough. Then the
    family $(V^{\varepsilon},
    \eta^{\varepsilon})_{0<\varepsilon < \varepsilon_0}$ is consistent with
    the system $(\mathcal{T}_b)$, where $V^{\varepsilon} = \nabla \psi^{\varepsilon}$.
    \end{prpstn}

    \vspace{1em}

    \noindent \begin{proof}
    First remark that since the velocity field $V^{\varepsilon}$ is
    irrotationnal, we have $\nabla \times \tilde{V}^{\varepsilon} = O(\varepsilon)$.
    And since $(\nabla\psi^{\varepsilon}, \eta^{\varepsilon})_{0
      < \varepsilon < \varepsilon_0}$ is consistent with the
    Boussinesq-like system $(\mathcal{B}_2)$, the previous
    computations yield directly the result. \end{proof}

    \vspace{1em} 

    \subsection{Derivation of a class of equivalent systems}

    In the previous section, we saw that a suitable change of
    variable comes from considering $V_{\theta}$, the horizontal
    component of the velocity at the height $-1+\theta$ ($\theta \in [0,1]$),
	instead of the horizontal component of the
    velocity field at the free surface. We can remark that the link between these two
    variables (and hence the adequate change of variables) can be
    derived from the expression of $u_{app}$ computed during the
    asymptotic expansion process of the operator
    $Z_{\varepsilon}(\varepsilon\eta,\beta b)$,
    which implies that we must adapt our change of
    variable for our regime of strong variations since the expression of $u_{app}$
    relies on the considered regime of bottom topography. \\

       \noindent Indeed, we
       saw in the previous section that the computation of the
       asymptotic developpment of $Z_{\varepsilon}(\eta,b)\psi$ relies
       on finding an approximate solution of the elliptic problem
       $(H)$ on the band $\mathcal{S} = [-1,0]\times\xR^2$. Starting
       from the truncation of the computed value of $u_{app}$ at the
       order $O(\varepsilon^2)$,
       $$
        u_{app} = \psi + \varepsilon \left[ \left(
        \frac{1}{2}-\frac{(z+1)^2}{2}\right) h^2\,\Delta \psi - z h
    \nabla h \cdot \nabla \psi \right] + O(\varepsilon^2) \;\;,
       $$
       where $\psi$ is the value of the velocity potential at the free
       surface, shows that $\nabla u_{app}(\cdot,z)$ gives an approximation
       at order $\varepsilon^2$ of the horizontal component of the
       velocity field, namely $V(\cdot,z) = \nabla \phi(\cdot,z)$ at height $z \in [-1,0]$.\\
       Consequently, in presence of huge bottom variations, the
    adequate change of variables is given by :
    $$
    V_{\theta} = \left[ 1-\frac{\varepsilon}{2} (\theta-1) (\theta
  \nabla(h^2\nabla \cdot\;) + \nabla\nabla \cdot (h^2\,.\;))\right]\,V\;\;,
    $$
    so that
    $$
    V = \left[ 1+\frac{\varepsilon}{2} (\theta-1) (\theta
  \nabla(h^2\nabla \cdot\;) + \nabla\nabla \cdot
  (h^2\,.\;))\right]\,V_{\theta} + O(\varepsilon^2)\;\;.
    $$
    From this change of variables, we easily compute the expressions of $\partial_t V$ and $\nabla \cdot \sqrt{h} V$ 
	which we plug into the system $(\mathcal{T}_b)$. By rewriting carefully the bottom terms in order to
	make the quantity $\sqrt{h}$ appear, one gets the following system :
	$$
        \left\{
         \begin{array}{l}
           \vspace{0.5em}
           \partial_t V_{\theta} +\sqrt{h}\;\nabla \eta  + \frac{\varepsilon}{2}\,
           \Big[\,F_h\left(\begin{array}{cc} V_\theta \\ \eta \end{array}\right) + \nabla\Big((\theta^2-1)
			h^2\nabla\cdot\partial_t V_{\theta}
	+ 2(\theta-1)h\nabla h \cdot \partial_t V_{\theta}\Big)\displaystyle\Big]= O(\varepsilon^2)\;\;,\\
           \vspace{0.5em}
           \partial_t \eta + \nabla (\sqrt{h}\cdot V_{\theta}) +
           \frac{\varepsilon}{2} \Big[\,f_h\left(\begin{array}{cc} V_\theta \\ \eta \end{array}\right) 
	-\nabla\cdot\left((\theta^2-\frac{1}{3})h^2\nabla\nabla\cdot(\sqrt{h}V_{\theta}) 
          +(\frac{3}{2}\,\theta^2-\frac{7}{6})h\nabla h\nabla\cdot(\sqrt{h}V_{\theta}) \right.\\
	\hspace{10em}\left.
	- \frac{(\theta-2)^2}{2}\sqrt{h}\nabla h(\nabla h \cdot V_{\theta}) - (\frac{\theta^2}{2}-2\theta + \frac{7}{6})
	h\sqrt{h}\nabla(\nabla h \cdot V_{\theta})
	 \right)\Big] = O(\varepsilon^2)\;\;.
         \end{array}
        \right.
        $$
	At this point a new problem arises. Applying the BBM trick in the exact same way as in
	the previous regime leads to a system that is never symmetric for any values of the parameters 
	$\theta, \lambda$ and $\mu$. Indeed, it implies to solve a numerical system
	on the unknowns $\theta$, $\lambda$ and $\mu$ which is
        over-determined. To deal with this problem, we simply introduce an
	additionnal unknown during the BBM trick process : we remark that the term $\partial_t V_{\theta}$ appears twice in the 
	dispersive terms of the first equation, we can then use two different expressions of $\partial_t V_{\theta}$, each with
	a different unknowns. This process is summed up in the following relations where we introduce the parameters $\lambda_1,
	\lambda_2$ and $\mu$ :
	$$
	\left\{
	\begin{array}{l}
	\vspace{1em}
	\partial_t V_{\theta} = (1-\lambda_1)\partial_t V_{\theta} - \lambda_1\sqrt{h}\nabla\eta + O(\varepsilon)\;\;,\\
	\vspace{1em}
	\partial_t V_{\theta} = (1-\lambda_2)\partial_t V_{\theta} - \lambda_2\sqrt{h}\nabla\eta + O(\varepsilon)\;\;,\\
	\nabla\cdot(\sqrt{h} V_{\theta}) = \mu \nabla \cdot(\sqrt{h} V_{\theta}) - (1-\mu)\partial_t \eta + O(\varepsilon)\;\;,
	\end{array}
	\right.
	$$
	where we use the first relation on the term $(\theta^2-1)h^2\nabla\cdot\partial_t V_{\theta}$ 
	and the second relation on the term $2(\theta-1)h\nabla h \cdot \partial_t V_{\theta}$.
	The last difficulty stands
	in the possibility to lose the possibly symmetric structure during the derivation of the final
	system. The key point is to rely on the quantity $\sqrt{h}$ which is crucial to write a symetric form of the dispersive
	terms. \\Finally, we formally derive a new class $(\mathcal{S}_b)$ of systems, 
	and we can prove that if a family $(V^{\varepsilon},
	\eta^{\varepsilon})_{0<\varepsilon<\varepsilon_0}$ is consistent with the system $\Gamma_b$ then 
	$(V^{\varepsilon}_{\theta}, 
	\eta^{\varepsilon})_{0<\varepsilon<\varepsilon_0}$ where $V_{\theta}^{\varepsilon} = \left[ 1-\frac{\varepsilon}{2} 
	(\theta-1) (\theta
	  \nabla(h^2\nabla \cdot\;) + \nabla\nabla \cdot (h^2\,.\;))\right]V^{\varepsilon}$ 
	is consistent with any of the following systems $(S_{\theta,\lambda_1,\lambda_2,\mu})$ :
	
	\vspace{1em}
	
	$$
        (S_{\theta,\lambda_1,\lambda_2,\mu})\left\{
         \begin{array}{l}
           \vspace{1em}
           \left(1-\frac{\varepsilon}{2} \mathcal{P}_{h}^1\right)\partial_t V +\sqrt{h}\;\nabla \eta  + \frac{\varepsilon}{2}\,
           \Big[\,F_h\left(\begin{array}{cc} V \\ \eta \end{array}\right) + b_1 \sqrt{h} \nabla \nabla \cdot(h^2
           \nabla\eta) + b_2 \sqrt{h} \nabla (h\nabla h \cdot \nabla
           \eta) \\
			\hspace{19em}\vspace{1em}  + b_3\nabla h \nabla 
			\cdot (h\sqrt{h}\nabla\eta) +b_4 \sqrt{h} \nabla h ( \nabla h \cdot \nabla \eta) \Big]= 0\;\;,\\
           \vspace{1em}
           \left(1-\frac{\varepsilon}{2} \mathcal{P}_{h}^2\right)\partial_t \eta + \nabla (\sqrt{h}\cdot V) +
           \frac{\varepsilon}{2} \Big[\,f_h\left(\begin{array}{cc} V \\ \eta \end{array}\right) 
	+c_1\nabla\cdot\Big(h^2\nabla\nabla\cdot(\sqrt{h}V) +
        c_2\nabla\cdot(h\nabla h\nabla\cdot(\sqrt{h}V)) \\
	\hspace{15em}\vspace{1em}  
	+ c_3\nabla\cdot(
	h\sqrt{h}\nabla(\nabla h \cdot V)) + c_4\nabla\cdot(\sqrt{h}\nabla h(\nabla h \cdot V))\Big)
	 \Big] = 0\;\;.
        \end{array}
        \right.
        $$
	where the operators $\mathcal{P}_{h}^1$ and $\mathcal{P}_{h}^2$ are defined by
	$$
	\left\{
	\begin{array}{l}
	\vspace{0.5em}
	\mathcal{P}_{h}^1 = (1-\theta)\Big((1-\lambda_1)(\theta+1)\nabla(h^2\nabla\cdot\hspace{1em})+2(1-\lambda_2)\nabla(
	h\nabla h\cdot\hspace{1em})\Big)\;\;,\\
	\mathcal{P}_{h}^2 = (1-\mu)\Big((\theta^2-\frac{1}{3})\nabla\cdot(h^2\nabla\hspace{1em})
	+(\frac{3}{2}\,\theta^2-\frac{7}{6})\nabla\cdot(
	h\nabla h\times\hspace{1em})\Big)\;\;,\\
	\end{array}
	\right.
	$$
	and the parameters $(a_i)_{1\le i \le 4}, (b_i)_{1\le i \le 4}$ have the following expressions :
	$$
	\left\{
	\begin{array}{ll}
	\vspace{0.5em}
	b_1= \lambda_1 (1-\theta^2); & c_1= \mu (\theta^2-\frac{1}{3}); \\
	\vspace{0.5em}
	b_2= (1-\theta)(2\lambda_2-\frac{3}{2}\lambda_1(1+\theta)); & c_2= \mu(\frac{3}{2}\theta^2-\frac{7}{6});\\
	\vspace{0.5em}
	b_3= \frac{\lambda_1}{2}(1-\theta^2); & c_3= -\frac{1}{2}\theta^2 + 2\theta-\frac{7}{6};\\
	b_4= (1-\theta)(\lambda_2-\frac{\lambda_1}{2}(1+\theta)); & c_4=\frac{1}{2}(\theta-2)^2; 
	\end{array}
	\right.
	$$
	
        \vspace{1em}

        \noindent The previous computations are summed up in the following
        proposition.

        \begin{prpstn}
          Let $\theta \in [0,1]$ and $(\psi^{\varepsilon}, \eta^{\varepsilon})_{0
            < \varepsilon < \varepsilon_0}$ a family of solutions of (\ref{S0}) such that $(\nabla\psi^{\varepsilon},
          \eta^{\varepsilon})_{0<\varepsilon<\varepsilon_0}$ is bounded with respect to $\varepsilon$ in
          $W^{1,\infty}([0,\frac{T}{\varepsilon}];\,H^{\sigma}(\xR^d)^{d+1})$ with $\sigma$ large enough. 
         We define $V^{\varepsilon} = \nabla \psi^{\varepsilon}$ and
          $\tilde{V}^{\varepsilon} =  \left(1-\frac{\varepsilon}{2} 
	(\theta-1) (\theta
	  \nabla(h^2\nabla \cdot\;) + \nabla\nabla \cdot
          (h^2\,.\;))\right)\,\left(\sqrt{h} +
         \frac{\varepsilon}{2}\,\frac{\eta}{\sqrt{h}} \right)\,V^{\varepsilon}$.
         Then for all $(\lambda_1,\lambda_2,\mu) \in \xR^3$, the
          family $(\tilde{V}^{\varepsilon},
          \eta^{\varepsilon})_{0<\varepsilon < \varepsilon_0}$ is consistent with
          the system $(S_{\theta,\lambda_1,\lambda_2,\mu})$.
        \end{prpstn}

        \vspace{1em}

        \noindent Moreover, we have the following proposition on the existence
        of a subclass of $(\mathcal{S}_b)$ composed with fully
        symmetric systems.

	\begin{prpstn} \label{existsym}
	There exists at least one value of $(\theta,\lambda_1,\lambda_2,\mu)$ 
	such that the system 
	$(S_{\theta,\lambda_1,\lambda_2,\mu})$ is fully symmetric.
	\end{prpstn}
	
	\vspace{1em}
	
	\noindent \begin{proof}
	We are concerned here with the resolution of the following system :
	$$
	\left\{
	\begin{array}{l}
	b_1 = c_1\;,\\
	b_2 = -c_2\;,\\
	b_3 = c_3\;,\\
	b_4 = -c_4\;.\\
	\end{array}
	\right.
	$$
	This system on $(\theta,\lambda_1,\lambda_2,\mu)$
        have at least one
	solution that gives the following approximate values :
	$$
	\left\{
	\begin{array}{l}
	\theta \approx 0.6318\;,\\
	\lambda_1 \approx -0.3416\;,\\
	\lambda_2 \approx -2.8209\;,\\
	\mu \approx -3.1157\;,\\
	\end{array}
	\right.
	$$
	which ends the proof. \end{proof}

	\vspace{1em}

        \noindent From now on, we only consider this solution and its
        approximate values.

        \vspace{1em}
	
	\subsection{The fully symmetric systems}
	
	Thanks to Proposition \ref{existsym}, we know that some of the 
	systems $(S_{\theta,\lambda_1,\lambda_2,\mu})$ 
	of the class $(\Lambda_b)$ are
	completely symmetric : we hence denote by $\Sigma_b$ the non-empty subclass of 
	$(\Lambda_b)$ composed with these
	symmetric systems. Unfortunately, we do not have the same kind of existence theory 
	on these systems as in the
	previous regime. \\
	Indeed, the main difference consists in the order one terms of the two equations
	$\left(\begin{array}{c} \sqrt{h}\nabla\eta \\ \nabla\cdot(\sqrt{h}V) \end{array} \right)
	$. In order to focus
	on the problem, we rewrite these terms :  
	$A(X,\partial_X)\left(\begin{array}{c} V \\ \eta 
	\end{array} \right)$ where $A(X,\partial_X) = \left(\begin{array}{cc} 0 & 
	\sqrt{h}\nabla_X \\  
	\nabla_X\cdot(\sqrt{h}\times\;) & 0 \end{array}\right)$. 
	The proof of the existence of solutions
	on a short time scale is not modified by these terms, the classical proof is still 
	valid. However, the fact that 
	the matrix $A$ depends on the bottom term $h$ is
	a real problem as far as the long time existence is concerned : indeed, one crucial 
	point of the proof here relies
	on the size of the quantity $\frac{\nabla h}{\varepsilon}$ on which we have no piece 
	of information. 
	The only case wherein we are surely able to demonstrate the long time existence is the 
	case where $\nabla h$ is of order
	$O(\varepsilon)$ : the term $\frac{\nabla h}{\varepsilon}$ is then of order $O(1)$ and 
	we can conclude. In all 
	other cases, the classical proof fails to provide a rigourous demonstration of the long 
	time existence of solutions
	to these symmetric systems. 
	Nevertheless, we are able to enounce the following proposition :
	
	\begin{prpstn}\label{symstrong}
        Let $s > \frac{d}{2}+1$ and $(\theta, \lambda_1, \lambda_2, \mu)$ be such that
         the system $(S_{\theta,\lambda_1, \lambda_2, \mu})$ belongs to the
         class $\Sigma_b$.\\
       Then for all $(V_0, \eta_0) \in H^s(\xR^d)^{d+1}$, there exists a
       time $T_0$ independant of $\varepsilon$ and a unique solution
       $(V,\eta) \in C([0,T_0]; H^s(\xR^d)^{d+1}) \cap\,C^1([0,T_0]; H^{s-3}(\xR^d)^{d+1}) $ to
       the system $(S_{\theta,\lambda_1, \lambda_2, \mu})$ such that
       $(V,\eta)_{|_{t=0}} = (V_0, \eta_0)$. \\
       Furthermore, this unique solution is bounded independently of
       $\varepsilon$ in the following sense : there exists a constant
       $C_0$ independent of $\varepsilon$ such that for all $k$
       verifying $s-3k > \frac{d}{2}+1$, we have :
       $$
       |(V,\eta)|_{W^{k,\infty}([0,T_0];H^{s-3k}(\xR^d)^{d+1})} \le C_0\;\;.
       $$
       Besides, if we suppose that $\nabla h = O(\varepsilon)$, the previous result becomes 
       valid on the long time
       interval $[0,\frac{T_0}{\varepsilon}]$.
	\end{prpstn}
	
	\vspace{1em}
	
	\noindent \begin{proof}
	The key point of the proof is to demonstrate that 
	the elliptic operator $1-\frac{\varepsilon}{2}
	\left(\begin{array}{cc} \mathcal{P}_h^1 \\ \mathcal{P}_h^2 \end{array}\right)$ 
	is a positive one. We first focus on
	$\mathcal{P}_h^1$ :
	$$
	(1-\frac{\varepsilon}{2}\mathcal{P}_h^1 V, V) = |V|_2^2+\frac{\varepsilon}{2}
	(1-\theta^2)(1-\lambda_1)|h \nabla
	\cdot V|_2^2 +\varepsilon (1-\theta) (1-\lambda_2)(\nabla h \cdot V, h \nabla \cdot V)
	$$
        Using the following inequality (satisfied for all $a \in \xR$)
        :
        $$
        \Big|(\nabla h \cdot V, h \nabla \cdot V)\Big| \le \frac{a^2}{2}
        |h\nabla\cdot V|_2^2 + \frac{1}{2a^2} |\nabla h \cdot V|_2^2\;\;,
        $$
        and taking $a^2=
        \frac{(1+\theta)(1-\lambda_1)}{1-\lambda_2}$
        leads to :
        $$
        (1-\frac{\varepsilon}{2}\mathcal{P}_h^1 V, V) \ge |V|_2^2 -
        \frac{\varepsilon}{2}\frac{(1-\theta)(1-\lambda_2)^2}{(1+\theta)(1-\lambda_1)} |\nabla h \cdot V |_2^2
        $$
        Using the classical Cauchy-Schwartz inequality leads finally
        to :
        $$
         (1-\frac{\varepsilon}{2}\mathcal{P}_h^1 V, V) \ge
         \Big(1-\frac{\varepsilon}{2}\frac{(1-\theta)(1-\lambda_2)^2}{(1+\theta)(1-\lambda_1)}|\nabla h|_2^2\Big)\,|V|_2^2\;\;,
        $$
        At this point, if we take a small enough value of
        $\varepsilon$, f.e. $\varepsilon \le \frac{2(1+\theta)(1-\lambda_1)}{(1-\theta)(1-\lambda_2)^2\,|\nabla
          h|_2^2}$, it ensures the global positivity of $\mathcal{P}_h^1$.
	On $\mathcal{P}_h^2$, we use the same method :
	$$
	(1-\frac{\varepsilon}{2}\mathcal{P}_h^2 \eta, \eta) = 
	|\eta|_2^2+\frac{\varepsilon}{2} (1-\mu) (\theta^2-\frac{1}{3})|h \nabla
	\eta|_2^2 +\frac{\varepsilon}{2} (1-\mu) (\frac{3}{2}\theta^2
        - \frac{7}{6}) (\eta \nabla h , 
	h \nabla \eta)
	$$
        Using the same ideas as previously, one gets :
        $$
         (1-\frac{\varepsilon}{2}\mathcal{P}_h^2 \eta, \eta) \ge
         \Big(1-\frac{\varepsilon}{8}\frac{(1-\mu)(\frac{3}{2}\theta^2-\frac{7}{6})^2}{\theta^2-\frac{1}{3}}\,|\nabla h|_2^2\Big)\,|\eta|_2^2\;\;,
        $$
        Once more, if we take f.e. $\varepsilon \le
        \frac{8(\theta^2-\frac{1}{3})}{(1-\mu)(\frac{3}{2}\theta^2-\frac{7}{6})^2\,|\nabla h|_2^2}$, we have the global positivity
        of $\mathcal{P}_h^2$.\\
        Consequently, taking $\varepsilon \le
        \min(\frac{2(1+\theta)(1-\lambda_1)}{(1-\theta)(1-\lambda_2)^2\,|\nabla
          h|_2^2},\frac{8(\theta^2-\frac{1}{3})}{(1-\mu)(\frac{3}{2}\theta^2-\frac{7}{6})^2\,|\nabla h|_2^2})$ ensures that the operator $1-\frac{\varepsilon}{2}
	\left(\begin{array}{cc} \mathcal{P}_h^1 \\ \mathcal{P}_h^2
          \end{array}\right)$ is positive.\\
        At this point, 
	using this result and performing usual energy
	estimates on the system proves the existence of a time $T$ such that there exists 
	an unique solution 
	$(V,\eta) \in C([0,T]; 
	H^s(\xR^d)^{d+1})\cap C^1([0,T]; H^{s-3}(\xR^d)^{d+1}) $
        to the system. \end{proof}
	
    \vspace{1em}
    
    \noindent This result gives us
    an efficient theoretical background to contruct approximate solutions of the water
    waves problem on a time scale $O(1)$, and $O(\frac{1}{\varepsilon})$ in the case $\nabla h = O(\varepsilon)$.\\
    This contruction follows the same steps - but in a different order - as the contruction
    of approximate solutions for the first regime : we consider a solution $(\psi^{\varepsilon},
    \eta^{\varepsilon})$ to the formulation (\ref{S0}) of the water
    waves problem. We take
    initial data $(\psi_0^{\varepsilon}, \eta_0^{\varepsilon})$ such
    that $(\nabla\psi_0^{\varepsilon}, \eta_0^{\varepsilon}) \in
    H^s(\xR^d)^{d+1}$ for a suitably large value of $s$. We then define
    $V^{\varepsilon} = \nabla \psi^{\varepsilon}$ and
    $V^{\varepsilon}_0 = \nabla \psi^{\varepsilon}_0$ : we first contruct the
    data $(V_{\Sigma,0}^{\varepsilon},
         \eta_{\Sigma,0}^{\varepsilon})$ by applying the two successive
    changes of variable on the data $(V_0^{\varepsilon}, \eta_0^{\varepsilon})$.
    We then choose the parameters $(\theta, \lambda_1, \lambda_2 \mu) \in
    [0,1]\times \xR^2$ such that the system
     $(T_{\theta,\lambda_1, \lambda_2, \mu}^{\;\;\;1})$ is completely symmetric. 
	Using Proposition \ref{symstrong},
    we know that there exists a unique solution to this system
    with initial data $(V_{\Sigma,0}^{\varepsilon}, \eta_{\Sigma,0}^{\varepsilon})$ : we denote this solution by
    $(V_{\Sigma}^{\varepsilon}, \eta_{\Sigma}^{\varepsilon})$.
    From this exact solution of the symmetric system
    $(S_{\theta,\lambda_1, \lambda_2, \mu})$, we finally construct an approximate solution of the water
    waves problem by successively and approximatively inverting the two changes of
    variable as shown below (which is possible if $\varepsilon$ is
    small enough) :
    
    \vspace{1em}
    
     $$ \label{approx2}
       \left\{
       \begin{array}{rcl}
         \vspace{0.5em}
       V_{app}^{\varepsilon} &=&  \left(\frac{1}{\sqrt{h}} -
         \frac{\varepsilon}{2}\,\frac{\eta^{\varepsilon}}{h\sqrt{h}}\right)
       \Big( 1+\frac{\varepsilon}{2} (\theta-1) (\theta
	  \nabla(h^2\nabla \cdot V_{\Sigma}^{\varepsilon}) + \nabla\nabla \cdot
	  (h^2\,V_{\Sigma}^{\varepsilon}))\Big)\\
       \eta_{app}^{\varepsilon} &=& \eta_{\Sigma}^{\varepsilon}
       \end{array}
       \right.      
    $$

       \vspace{1em}

       \noindent We are now able to enounce our final result :

       \vspace{1em}

       \begin{thrm}\label{strong}
         Let $T_1 \ge 0$, $s \ge \frac{d}{2}+1$, $\sigma \ge s+3$ and
         $(\nabla\psi_0^{\varepsilon},\eta_0^{\varepsilon})$ be in $H^{\sigma}(\xR^d)^{d+1}$.
         Let $(\psi^{\varepsilon},
         \eta^{\varepsilon})_{0 < \varepsilon < \varepsilon_0}$ be a
         family of
         solutions of (\ref{S0}) with initial data
         $(\nabla\psi_0^{\varepsilon},\eta_0^{\varepsilon})$ and such that $(\nabla\psi^{\varepsilon},
          \eta^{\varepsilon})_{0<\varepsilon<\varepsilon_0}$ is bounded with
         respect to $\varepsilon$ in
         $W^{1,\infty}([0,T_1]; H^{\sigma}(\xR^d)^{d+1})$.
         We define $V^{\varepsilon} = \nabla\psi^{\varepsilon}$ and choose
         $(\theta, \lambda_1, \lambda_2 \mu) \in [0,1]\times\xR^2$  such that
         the system $(S_{\theta,\lambda_1, \lambda_2, \mu}^{\;\;\;1}) \in \Sigma$.\\
         Then for all $\varepsilon < \varepsilon_0$, there exists a
         time $T \le T_1$ such that we have :
         $$
         |V^{\varepsilon} -
         V^{\varepsilon}_{app}|_{L^{\infty}([0,T];H^s)} +
         |\eta^{\varepsilon} -
         \eta^{\varepsilon}_{app}|_{L^{\infty}([0,T];H^s)} \le
         C\,\varepsilon^2
         $$
         
         \vspace{0.5em}
         
         \noindent Besides, if we suppose that $\nabla h = O(\varepsilon)$ 
		 then $(V_{app}^{\varepsilon}, \eta_{app}^{\varepsilon})$
         approximates the water waves solutions on a long time scale :
         $$
         \forall t \in [0,\frac{T}{\varepsilon}]\;,\;\;\;|V^{\varepsilon} -
         V^{\varepsilon}_{app}|_{L^{\infty}([0,t];H^s)} +
         |\eta^{\varepsilon} -
         \eta^{\varepsilon}_{app}|_{L^{\infty}([0,t];H^s)} \le
         C\,\varepsilon^2 t
         $$
       \end{thrm}
       
       \vspace{2em}

       \noindent \begin{proof}
       The proof is an adaptation of the one of Theorem \ref{small}, and we omit it here. \end{proof}
       
	\vspace{1em}
       
	\begin{rmrk} In the general case, where we have no piece of information 
	on the size of the quantity
		$\frac{\nabla h}{\varepsilon}$, our analysis is complete on a short time scale.
		We have indeed an approximation on this interval of time, and
		we know from Lannes \cite{Lannes}
		the existence of solutions to the water waves problem on a short time scale 
		in 2D and
		3D. However, in the case $\nabla h = O(\varepsilon)$, 
	    this analysis is only complete in 2D - like in the first regime - since
		we do not know about the existence of solutions to the water waves problem 
		on a long time scale. \end{rmrk}

        \vspace{3em}

	\begin{acknowledgement} This
        work was supported by the ACI Jeunes chercheurs du minist\`ere
        de la Recherche ``Dispersion et nonlin\'earit\'es''.
        \end{acknowledgement}

        \vspace{1em}


\end{document}